# EXTENDING THE SCOPE OF EMPIRICAL LIKELIHOOD

By Nils Lid Hjort, Ian W. McKeague[1] and Ingrid Van Keilegom[2]

*University of Oslo, Columbia University, and*
*Université catholique de Louvain and Tilburg University (CentER)*

This article extends the scope of empirical likelihood methodology in three directions: to allow for plug-in estimates of nuisance parameters in estimating equations, slower than $\sqrt{n}$-rates of convergence, and settings in which there are a relatively large number of estimating equations compared to the sample size. Calibrating empirical likelihood confidence regions with plug-in is sometimes intractable due to the complexity of the asymptotics, so we introduce a bootstrap approximation that can be used in such situations. We provide a range of examples from survival analysis and nonparametric statistics to illustrate the main results.

**1. Introduction.** Empirical likelihood [Owen (1990, 2001)] has traditionally been used for providing confidence regions for multivariate means and, more generally, for parameters in estimating equations, under various standard assumptions: the number of estimating equations is fixed, they do not involve nuisance parameters, and the parameters of interest are estimable at $\sqrt{n}$-rate, where $n$ is the sample size. Under such assumptions and with i.i.d. observations [or even dependent observations; see, e.g., Chapter 8 of Owen (2001)], empirical likelihood (EL) based confidence regions can be calibrated using a nonparametric version of Wilks's theorem involving a chi-squared limiting distribution.

The aim of the present paper is to develop adaptations when the traditional assumptions are violated. More specifically, under certain asymptotic

Received December 2006; revised August 2007.

[1]Supported by NSF Grant DMS-05-05201.

[2]Supported by IAP Research Network Grants P5/24 and P6/03 of the Belgian government.

stability conditions, we establish generalizations of the basic theorem of EL to allow for plug-in estimates of nuisance parameters in the estimating equations, for slower than $\sqrt{n}$-rates of convergence, and for i.i.d. settings in which there are a relatively large number of estimating equations compared to the sample size. Several of our examples share the characteristic that they would be harder to analyze with other methods. In particular, the method of profile EL [see, e.g., Owen (2001), page 42] for dealing with nuisance parameters in estimating equations is often not applicable for infinite-dimensional nuisance parameters, and even when it is applicable, implementation can be computationally difficult. The triangular array EL theorem of Owen [(2001), page 85] applies under slower than $\sqrt{n}$-rates, and has been useful in the context of nonparametric density estimation, for instance, but is not flexible enough to handle estimating functions with plug-in.

The use of plug-in for nuisance parameters in EL confidence regions is not new. It has recently been applied in various survival analysis contexts; see Qin and Jing (2001a, 2001b), Wang and Jing (2001), Li and Wang (2003) and Qin and Tsao (2003). The technique has also been used in survey sampling with imputation for missing response; see Wang and Rao (2002). Our aim here, however, is to provide a more widely applicable version of this approach, that can accommodate a wide array of examples, allowing both plug-in and slower than $\sqrt{n}$-rates of convergence. We take the point of view that it is preferable to derive a general result using generic assumptions, that can be checked in a large number of applications, rather than reinventing the basic theory on each occasion. Calibrating EL confidence regions with plug-in is sometimes intractable due to the complexity of the asymptotics, so we introduce a bootstrap approximation that can be used in such situations.

To illustrate our general results we consider a range of examples from survival analysis and nonparametric statistics in settings where the inference is based on estimating functions. In particular, we look at functionals of survival distributions with right censored data [treated via EL in Wang and Jing (2001)], the error distribution in nonparametric regression [Akritas and Van Keilegom (2001)], density estimation [treated by EL in Hall and Owen (1993) and Chen (1996)], and survival function estimation from current status data [van der Vaart and van der Laan (2006)].

Standard maximum likelihood theory for parametric models, as well as EL theory, keeps the dimension of the parameter (or the number of estimating equations) fixed, say at $p$, as sample size $n$ grows. This is what leads to asymptotic normality, Wilks type theorems for likelihood ratio statistics and Owen type theorems for EL. Portnoy (1986, 1988) and others have investigated the extent to which maximum likelihood theory based results still hold, when $p$ is allowed to increase with $n$. The canonical growth restriction for normal approximations to hold is that $p^2/n \to 0$, while $p^{3/2}/n \to 0$ typically suffices for certain quadratic approximations associated with Wilks theorems to hold.

In this article we investigate the similar problem of finding conditions under which the EL methods continue to work adequately when $p$ grows. The canonical growth condition will be seen to be $p^3/n \to 0$. Under this condition, in addition to other requirements that have to do with stability of eigenvalues of covariance matrices, minus twice the log-EL can be approximated well enough with a certain quadratic form that in itself is close to a $\chi_p^2$.

We should add that in situations with a high number of parameters the typical aim is not to provide a simultaneous confidence region for the full parameter vector, say $(\mu_1, \ldots, \mu_p)$. It could rather be to test whether a subset of the parameters have zero values, or to compare one distribution with another, or, more generally, to make inference for a focus parameter of dimension $q < p$, say $f(\mu_1, \ldots, \mu_p)$. For any linear map $f$, these tasks can be carried out inside our framework for growing $p$ by constructing a $q$-dimensional confidence region in which $q$ grows with $n$. For further discussion in the context of a regression example, see Section 5.4.

The paper is organized as follows. Section 2 develops the EL theory with plug-in and the bootstrap approximation of the limiting distribution of the EL statistic. Six examples, including two involving slower than $\sqrt{n}$-rates of convergence, are discussed in Section 3. In Section 4 we examine the limiting behavior of the EL statistic in situations where the number of estimating functions is allowed to increase with growing sample size. Some examples are presented in Section 5, including setups with "growing polynomial regression" and "growing exponential families." Proofs can be found in the Appendix.

## 2. Plug-in empirical likelihood.

We first describe the general framework. The basic idea of empirical likelihood (EL) is to regard the observations $X_1, \ldots, X_n$ as if they are i.i.d. from a fixed and unknown $d$-dimensional distribution $P$, and to model $P$ by a multinomial distribution concentrated on the observations. Inference for the parameter(s) of interest, $\theta_0 = \theta(P) \in \Theta$, is then carried out using a $p$-dimensional estimating function of the form $m_n(X, \theta, h)$, where, for the purposes of the present paper, $h$ is a (possibly infinite-dimensional) "nuisance" parameter with unknown true value $h_0 = h(P) \in \mathcal{H}$.

When $h_0$ is known, it can replace $h$ in the EL ratio function

$$\text{EL}_n(\theta, h) = \max\left\{\prod_{i=1}^n (nw_i): \text{ each } w_i > 0, \sum_{i=1}^n w_i = 1, \sum_{i=1}^n w_i m_n(X_i, \theta, h) = 0\right\},$$

leading to a confidence region $\{\theta: \text{EL}_n(\theta, h_0) > c\}$ for $\theta_0$, where $c$ is a suitable positive constant, and the maximum of the empty set is defined to be zero. The constant $c$ can be calibrated using Owen's (1990) EL theorem, provided $m_n = m$ does not depend on $n$: if the observations are i.i.d. and $m(X, \theta_0, h_0)$ has zero mean and a positive definite covariance matrix, then $-2 \log \mathrm{EL}_n(\theta_0, h_0) \to_d \chi^2_p$, where $\chi^2_p$ has a chi-squared distribution with $p$ degrees of freedom.

2.1. *Main result.* We now establish a generalization of Owen's result in which the unknown $h_0$ is replaced by an estimator $\hat{h}$, and the estimating function is allowed to depend on $n$. This result will provide a way of calibrating $\{\theta : \mathrm{EL}_n(\theta, \hat{h}) > c\}$ as a confidence region for $\theta_0$. We extract the basic structure of Owen's result, and only impose an existence condition, (A0) below, and some "generic" asymptotic stability conditions, (A1)–(A3) below. These conditions ensure a nondegenerate limiting distribution, but do not require i.i.d. observations or consistency of $\hat{h}$, although such structure may very well be helpful for checking the conditions in specific applications. Our proof (placed in the Appendix) uses tools somewhat different from those usually employed in the EL literature, as in, for example, Owen (2001), Chapter 11; see also Remark 2.7 below.

We use the following notation throughout. For vectors $v$, let $\|v\|$ denote the Euclidean norm, and $v^{\otimes 2} = vv^{\mathrm{t}}$. For matrices $V = (v_{i,j})$, let $|V| = \max_{i,j} |v_{i,j}|$.

Let $\{a_n\}$ be a sequence of positive constants bounded away from zero, and $U$ a nondegenerate $p$-dimensional random vector. In most of the applications we consider, $a_n = 1$ and $U \sim \mathrm{N}_p(0, V_1)$, where the covariance matrix $V_1$ is positive definite, but the extra generality can be useful in some applications. Let $V_2$ denote a $p \times p$ positive definite covariance matrix. The following conditions are needed:

(A0)  $P\{\mathrm{EL}_n(\theta_0, \hat{h}) = 0\} \to 0$.

(A1)  $\sum_{i=1}^n m_n(X_i, \theta_0, \hat{h}) \to_d U$.

(A2)  $a_n \sum_{i=1}^n m_n^{\otimes 2}(X_i, \theta_0, \hat{h}) \to_{\mathrm{pr}} V_2$.

(A3)  $a_n \max_{1 \le i \le n} \|m_n(X_i, \theta_0, \hat{h})\| \to_{\mathrm{pr}} 0$.

As pointed out by a referee, $\hat{h}$ just plays the role of indicating that $m_n$ is being estimated, and we could replace $m_n(X, \theta, \hat{h})$ by the simpler notation $\hat{m}_n(X, \theta)$. This also covers situations in which $h$ depends on $\theta$ with an estimating function of the form $m_n(X, \theta, \hat{h}_\theta)$. We prefer to include $\hat{h}$ explicitly in the notation, however, because all our examples involve a plug-in estimator, as does our bootstrap result in Section 2.3.

Condition (A0) is equivalent to $P(0 \in \mathcal{C}_n) \to 1$, where $\mathcal{C}_n$ denotes the interior of the convex hull of $\{m_n(X_i, \theta_0, \hat{h}), \ i = 1, \ldots, n\}$ and 0 is the zero vector in $\mathbb{R}^p$. This is the basic existence condition needed for EL to be useful in our general setting. Below we describe how the EL statistic can be expressed, up to a negligible remainder term, as a quadratic form involving the left-hand sides of (A1) and (A2), so these conditions play a natural role in the asymptotics (see Remark 2.7). Finally, (A3) is required to obtain the negligibility of the remainder term. For the practical verification of these conditions, we refer the reader to Section 3, where they are checked in detail in a number of applications.

THEOREM 2.1. *If* (A0)–(A3) *hold, then* $-2a_n^{-1} \log \mathrm{EL}_n(\theta_0, \widehat{h}) \to_d U^{\mathrm{t}} V_2^{-1} U$.

2.2. *Remarks.* This theorem is related to many results in the literature, which we now discuss, along with a sketch of its proof; the complete proof appears in the Appendix.

REMARK 2.1. Owen's EL theorem follows from Theorem 2.1 by taking $a_n = 1$ and $m_n = m/\sqrt{n}$. Indeed, (A0) then holds using an argument involving the Glivenko–Cantelli theorem over half-spaces [see page 219 of Owen (2001)], (A1) by the multivariate central limit theorem, (A2) by the law of large numbers, and (A3) by a Borel–Cantelli argument [Lemma 11.2 of Owen (2001)].

REMARK 2.2. When $U \sim \mathrm{N}_p(0, V_1)$ with $V_1$ positive definite, the limit distribution above may be expressed as $r_1 \chi_{1,1}^2 + \cdots + r_p \chi_{1,p}^2$, where the $\chi_{1,j}^2$'s are independent chi-squared random variables with one degree of freedom and the weights $r_1, \ldots, r_p$ are the eigenvalues of $V_2^{-1} V_1$; cf. Lemma 3 of Qin and Jing (2001a). If, in addition, $V_1$ and $V_2$ coincide, we have the standard $\chi_p^2$ limit distribution. When $V_1$ and $V_2$ are not identical, the weights $r_1, \ldots, r_p$ may need to be estimated, for example via consistent estimators $\widehat{V}_1$, $\widehat{V}_2$ and computing the eigenvalues of $\widehat{V}_2^{-1} \widehat{V}_1$. It is not possible to say anything in general about estimation of $V_1$, which will depend on the structure of the specific application; later in this section we examine a bootstrap approach which can be applied when $V_1$ is difficult to estimate by other means. For $a_n = 1$, an estimator of $V_2$ is easily provided given plug-in of a consistent estimator $\widehat{\theta}$ for $\theta_0$. In the Appendix we show that $\widehat{V}_2 = \sum_{i=1}^{n} m_n^{\otimes 2}(X_i, \widehat{\theta}, \widehat{h})$ consistently estimates $V_2$ under the following two additional conditions: there exists a $p \times p$-matrix-valued function $V(\theta, h)$ such that

(A4) For some subset $\bar{\mathcal{H}}$ of $\mathcal{H}$ such that $P\{\widehat{h} \in \bar{\mathcal{H}}\} \to 1$, and for some $\delta > 0$,

$$\sup_{\|\theta - \theta_0\| < \delta, h \in \bar{\mathcal{H}}} \left| \sum_{i=1}^{n} m_n^{\otimes 2}(X_i, \theta, h) - V(\theta, h) \right| \to_{\mathrm{pr}} 0;$$

(A5) $\sup_{\|\theta - \theta_0\| \leq \delta_n, h \in \bar{\mathcal{H}}} |V(\theta, h) - V(\theta_0, h)| \to 0$ for any real sequence $\delta_n \downarrow 0$.

When the observations are i.i.d. and $m_n = m/\sqrt{n}$ for some function $m(X, \theta, h)$ that does not depend on $n$, we would expect to use $V(\theta, h) = \mathrm{E} m^{\otimes 2}(X_1, \theta, h)$ and then (A4) amounts to a (convergence-in-probability) version of the Glivenko–Cantelli property for $\mathcal{F} = \{m^{\otimes 2}(\cdot, \theta, h) \colon \|\theta - \theta_0\| < \delta, h \in \bar{\mathcal{H}}\}$.

REMARK 2.3.   For i.i.d. observations and $m_n = m/\sqrt{n}$, with $m(X, \theta_0, h_0)$ having zero mean and a finite covariance matrix $V_0$, the multivariate central limit theorem implies that $\sum_{i=1}^n m_n(X_i, \theta_0, h_0)$ tends to $\mathrm{N}_p(0, V_0)$, so condition (A1) describes the perturbation of $V_0$ due to replacing $h_0$ by $\hat{h}$. In the "highly smooth" case that $M(\theta_0, \hat{h}) = o_{\mathrm{pr}}(n^{-1/2})$, where $M(\theta, h) = \mathrm{E} m(X, \theta, h)$, it can be shown (under some additional assumptions) that there is no perturbation: $V_1 = V_0$. For instance, suppose that the class of functions $\{m(\cdot, \theta_0, h) \colon h \in \mathcal{H}\}$ is Donsker, and $\hat{h}$ is consistent in the sense that $\rho_j(\hat{h}, h_0) \to_{\mathrm{pr}} 0$ for $j = 1, \ldots, p$, where $\rho_j(h, h_0)$ is the $L^2(P)$ distance between $m_j(X, \theta_0, h)$ and $m_j(X, \theta_0, h_0)$. Then

$$\sum_{i=1}^n m_n(X_i, \theta_0, \hat{h}) = n^{-1/2} \sum_{i=1}^n \{m(X_i, \theta_0, \hat{h}) - M(\theta_0, \hat{h})\} + \sqrt{n} M(\theta_0, \hat{h})$$

tends to $\mathrm{N}_p(0, V_0)$, so $V_1 = V_0$, where empirical process theory is used to obtain weak convergence of the first term; cf. van der Vaart (1998), page 280. However, $M(\theta_0, \hat{h}) = o_{\mathrm{pr}}(n^{-1/2})$ is a strong condition, so we have avoided using it in favor of the less restrictive condition (A1), which is flexible enough to be checked within the context of the examples considered in the next section.

REMARK 2.4.   Kitamura (1997) introduces blockwise EL with estimating functions, without plug-in, in models having weakly dependent stationary observations. The maximum EL estimator under blocking is shown to have greater efficiency than the standard maximum EL estimator, but the blockwise approach has not been extended to allow plug-in. Standard EL (with plug-in), however, can still provide accurate confidence sets under dependent observations, for according to Theorem 2.1 the limiting distribution of the standard EL statistic, while not chi-square, is of a tractable form. If $m_n = m/\sqrt{n}$ and there is no plug-in, conditions (A1) and (A2) can be checked by central limit theorems and ergodic theorems for weakly dependent sequences. Condition (A3) holds provided $\mathrm{E} \|m(X, \theta_0)\|^2 < \infty$ by a Borel–Cantelli argument [cf. Owen (2001), Lemma 11.2]. For an estimating function $m(X, \theta)$ such that $\mathrm{E} m(X, \theta_0) = 0$, the limiting distribution of the EL statistic is as in Remark 2.2 with $V_1 = \sum_{i=1}^{\infty} \mathrm{Cov}\{m(X_1, \theta_0), m(X_i, \theta_0)\}$ and $V_2 = \mathrm{Var}\{m(X, \theta_0)\}$, which could be estimated easily.

REMARK 2.5. Nordman, Sibbertsen and Lahiri (2007) develop block-wise EL for the mean of the long-range dependent (stationary and ergodic) process $X_i = G(Z_i)$, where $\{Z_i\}$ is a stationary sequence of $N(0, 1)$ random variables such that $\text{cov}(Z_i, Z_{i+n}) = n^{-\alpha} L(n)$, for some $0 < \alpha < 1$ and slowly varying $L(\cdot)$, and $G(\cdot)$ is a Borel function with $G(Z_1)$ having finite mean $\theta_0$ and finite variance $\sigma^2$. Suppose that $\alpha$, $L(\cdot)$ and $G(\cdot) - \theta_0$ are known and we use an estimating function of the form $m_n(X_i, \theta) = b_n(X_i - \theta)$, where $b_n$ depends on the rate of convergence of the sample mean of the $X_i$. Condition (A1) is checked using a result of Taqqu (1975), which shows that $b_n \sum_{i=1}^n (X_i - \theta_0) \to_d U$ if we specify $b_n = n^{\alpha/2-1} L(n)^{-1/2}$. Here $U$ is defined by a multiple Wiener integral and does not depend on $\theta_0$. Condition (A2) is checked by setting $a_n = n^{-1} b_n^{-2} = n^{1-\alpha} L(n)$ and using the ergodic theorem:

$$a_n \sum_{i=1}^n m_n(X_i, \theta_0)^2 = n^{-1} \sum_{i=1}^n (X_i - \theta_0)^2 \to_{\text{a.s.}} \sigma^2 = V_2.$$

In this case the choice of $a_n$ tends to infinity, and it is not possible to arrange $a_n = 1$.

REMARK 2.6. In the special case that the nuisance parameter $h$ is finite dimensional, the profile EL statistic

$$-2 \log \left\{ \max_h \text{EL}_n(\theta_0, h) \Big/ \max_{\theta, h} \text{EL}_n(\theta, h) \right\} \to_d \chi_q^2$$

under various regularity conditions [Qin and Lawless (1994), Corollary 5], where $q$ is the dimension of $\theta$. This provides an attractive method of obtaining an EL confidence region for $\theta$, and is easier than using plug-in, but it is restricted to finite-dimensional nuisance parameters and the estimating function needs to be differentiable in $(\theta, h)$. Bertail (2006) extended this approach to infinite-dimensional $h$ in some "highly smooth" cases (cf. Remark 2.3).

REMARK 2.7. Our proof of Theorem 2.1 differs from the usual EL approach in that we take the dual problem perspective; see, for example, Christianini and Shawe-Taylor (2000), Section 5.2, for the relevant convex optimization theory. An outline of the proof is as follows. Write $X_{n,i} = m_n(X_i, \theta_0, \hat{h})$. By (A0), with probability tending to 1, $\text{EL}_n = \text{EL}_n(\theta_0, \hat{h}) = \prod_{i=1}^n (1 + \hat{\lambda}^t X_{n,i})^{-1}$, where the $p$-vector of Lagrange multipliers $\hat{\lambda}$ satisfies $\sum_{i=1}^n X_{n,i}/(1 + \hat{\lambda}^t X_{n,i}) = 0$, as in Owen (2001), page 219. Thus, with probability tending to 1, we can express the EL statistic in dual form as

$$(1) \qquad -2 \log \text{EL}_n = G_n(\hat{\lambda}) = \sup_\lambda G_n(\lambda),$$

where $G_n(\lambda) = 2\sum_{i=1}^n \log(1 + \lambda^{\mathrm{t}} X_{n,i})$, and the domain of $G_n$ is the set on which it is defined (regarding $\log x$ as undefined for $x \leq 0$). Note here that $G_n$ is concave and achieves its maximum at $\widehat{\lambda}$ since $\nabla G_n(\widehat{\lambda}) = 0$. Now consider the following quadratic approximation to $G_n$:

$$G_n^*(\lambda) = 2\lambda^{\mathrm{t}} U_n - \lambda^{\mathrm{t}} V_n \lambda \qquad \text{where } U_n = \sum_{i=1}^n X_{n,i}, V_n = \sum_{i=1}^n X_{n,i}^{\otimes 2},$$

and the domain of $G_n^*$ is taken as the whole of $\mathbb{R}^p$. We show in our Appendix that the difference between the maxima of $G_n$ and $G_n^*$ (over their respective domains) is of order $o_{\mathrm{pr}}(a_n)$. Thus, by (1) and the fact that $G_n^*$ is maximized at $\lambda^* = V_n^{-1} U_n$ when $V_n$ is invertible (which happens with probability tending to 1), it follows that

$$(2) \qquad -2a_n^{-1} \log \mathrm{EL}_n = a_n^{-1} \sup_\lambda G_n^*(\lambda) + o_{\mathrm{pr}}(1) = U_n^{\mathrm{t}}(a_n V_n)^{-1} U_n + o_{\mathrm{pr}}(1),$$

which tends in distribution to $U^{\mathrm{t}} V_2^{-1} U$, via assumptions (A1) and (A2). It also follows from the proof that Theorem 2.1 continues to hold in cases where $(U_n, V_n) \to_d (U, V_2)$, with a random rather than a fixed $V_2$.

2.3. *Bootstrap calibration.* As mentioned above, the estimation of $V_1$ can be difficult in certain situations and, more seriously, $U$ may not be normally distributed, in which case a bootstrap calibration is desirable. The procedure developed below consists in replacing $U$ by a bootstrap approximation, and in consistently estimating $V_2$.

We restrict attention to i.i.d. data and $m_n = m/\sqrt{n}$. Assume that $M(\theta, h_0) = 0$ if and only if $\theta = \theta_0$, where $M(\theta, h) = \mathrm{E}m(X, \theta, h)$, and denote $M_n(\theta, h) = n^{-1}\sum_{i=1}^n m(X_i, \theta, h)$. Let $\{X_1^*, \ldots, X_n^*\}$ be drawn randomly with replacement from $\{X_1, \ldots, X_n\}$, let $\widehat{h}^*$ be the same as $\widehat{h}$ but based on the bootstrap data, and define $M_n^*(\theta, h) = n^{-1}\sum_{i=1}^n m(X_i^*, \theta, h)$. Also, let $\widehat{\theta}$ be a consistent estimator of $\theta_0$, and $\widehat{V}_2 = n^{-1}\sum_{i=1}^n m^{\otimes 2}(X_i, \widehat{\theta}, \widehat{h})$.

We use the abbreviated notation $\Delta_n = M_n - M$, as a function of $(\theta, h)$, and $\Delta_n^*$ denotes the bootstrap version of $\Delta_n$ (here and in the sequel we define the bootstrap version of any statistic as the expression obtained by replacing $M, M_n, \theta_0, h_0$ and $\widehat{h}$ by $M_n, M_n^*, \widehat{\theta}, \widehat{h}$ and $\widehat{h}^*$, resp.). Let $\mathcal{H}$ be a vector space of functions endowed with a pseudo-metric $\|\cdot\|_{\mathcal{H}}$, which is a sup-norm metric with respect to the $\theta$-argument and a pseudo-metric with respect to all the other arguments. Also let $\Phi_n = \sqrt{n}\{\Delta_n(\theta_0, h_0) + \Gamma(\theta_0, h_0)[\widehat{h} - h_0]\}$, where $\Gamma(\theta_0, h_0)[\widehat{h} - h_0]$ is the Gâteaux derivative of $M(\theta_0, h_0)$ in the direction $\widehat{h} - h_0$ [see, e.g., Bickel, Klaassen, Ritov and Wellner (1993), page 453]. The bootstrap analogue of $\Phi_n$ is denoted by $\Phi_n^*$. Finally, let $P^*$ denote the bootstrap distribution conditional on the data. The following conditions are needed to formulate the validity of the bootstrap approximation:

(B1) $\sup_{t \in \mathbb{R}^p} |P^*\{\Phi_n^* \le t\} - P\{\Phi_n \le t\}| \to_{\mathrm{pr}} 0$.

(B2) $\sup_{\|\theta - \theta_0\| \le \delta_n, \|h - h_0\|_{\mathcal{H}} \le \delta_n} \|\Delta_n(\theta, h) - \Delta_n(\theta_0, h_0)\| = o_{\mathrm{pr}}(n^{-1/2})$ for all $\delta_n \downarrow 0$.

(B3) $\|M(\theta_0, \widehat{h}) - M(\theta_0, h_0) - \Gamma(\theta_0, h_0)[\widehat{h} - h_0]\| \le c\|\widehat{h} - h_0\|_{\mathcal{H}}^2$ for some $c > 0$.

(B4) $\|\widehat{h} - h_0\|_{\mathcal{H}} = o_{\mathrm{pr}}(n^{-1/4})$.

(B5) The bootstrap analogues of conditions (B2)–(B4) hold pr-a.s.

THEOREM 2.2. *Under conditions* (A0)–(A5) *and* (B1)–(B5),

$$\sup_{t \ge 0} |P^*\{n[M_n^*(\widehat{\theta}, \widehat{h}^*) - M_n(\widehat{\theta}, \widehat{h})]^{\mathrm{t}} \widehat{V}_2^{-1}[M_n^*(\widehat{\theta}, \widehat{h}^*) - M_n(\widehat{\theta}, \widehat{h})] \le t\}$$

$$- P\{-2 \log \mathrm{EL}_n(\theta_0, \widehat{h}) \le t\}| \to_{\mathrm{pr}} 0.$$

REMARK 2.8. When $\widehat{\theta}$ is defined as the minimizer of $\|M_n(\theta, \widehat{h})\|$, sufficient conditions for $\widehat{\theta}$ to be consistent can be found in Theorem 1 in Chen, Linton and Van Keilegom (2003). In order to verify condition (B2) in the case of i.i.d. observations, it suffices by Corollary 2.3.12 in van der Vaart and Wellner (1996) to show that the class $\{m(\cdot, \theta, h) : \theta \in \Theta, h \in \mathcal{H}\}$ is Donsker, and that

$$\mathrm{Var}\{m(X, \theta, h) - m(X, \theta_0, h_0)\} \le K_1\|\theta - \theta_0\| + K_2\|h - h_0\|_{\mathcal{H}} + \varepsilon_n$$

for some $K_1, K_2 \ge 0$, and for some $\varepsilon_n \downarrow 0$. The former condition can be verified by making use of Theorem 3 in Chen, Linton and Van Keilegom (2003). The bootstrap analogue of (B2) then follows from Giné and Zinn (1990), provided

$$\mathrm{Var}^*\{m(X^*, \theta, h) - m(X^*, \widehat{\theta}, \widehat{h})\} \le K_1'\|\theta - \widehat{\theta}\| + K_2'\|h - \widehat{h}\|_{\mathcal{H}} + \varepsilon_n'$$

for some $K_1', K_2' = O(1)$ a.s. and for some $\varepsilon_n' = o(1)$ a.s. Finally, condition (B3) and its bootstrap version can often be verified by using a two-term Taylor expansion of $M(\theta_0, \widehat{h})$ and of $M(\widehat{\theta}, \widehat{h}^*)$ around $h_0$ and $\widehat{h}$, respectively.

**3. Applications of the plug-in theory.** This section gives six illustrations of the preceding plug-in theory. The first uses parametric plug-in for a nonparametric estimand while the five others effectively use nonparametric plug-in to solve nonparametric empirical likelihood problems. The last two are examples of situations where the rate of convergence of the estimator of $\theta_0$ is slower than the usual root-$n$ rate. All the examples use $a_n = 1$.

3.1. *Symmetric distribution functions.* Let $F$ be a continuous distribution function of a random variable $X$, that is symmetric about an unknown location $a$, so $F(x) = 1 - F(2a - x)$ for all $x$. Consider estimation of $\theta_0 = F(x)$ at a fixed $x$ from $n$ i.i.d. observations from $F$. The estimating function has $p = 2$ components (the first being the usual estimating function and the second making use of the symmetry assumption): $m_n = n^{-1/2}m$, with

$$m(X, \theta, a) = \begin{pmatrix} I\{X \leq x\} - \theta \\ I\{X > 2a - x\} - \theta \end{pmatrix}.$$

The plug-in estimator of $a$ is taken as the sample median $\hat{a}$. Let $\eta_0 = \min(\theta_0, 1 - \theta_0)$ and suppose $0 < \theta_0 < 1$. Condition (A2) holds and

$$V_2 = \begin{pmatrix} \theta_0(1 - \theta_0) & -\eta_0^2 \\ -\eta_0^2 & \theta_0(1 - \theta_0) \end{pmatrix}$$

when $\theta_0 \neq 1/2$, and $V_2$ is singular when $\theta_0 = 1/2$. A consistent estimator of $V_2$ is obtained by replacing $\theta_0$ by $\hat{F}(x)$, where $\hat{F}$ is the empirical distribution function of $X$. The validity of condition (A3) is straightforward. Now, let us turn to condition (A1). First note that

$$\sqrt{n}\{1 - \hat{F}(2\hat{a} - x) - \theta_0\}$$
$$= \sqrt{n}\{1 - F(2\hat{a} - x) - \hat{F}(2a - x) + F(2a - x) - \theta_0\} + o_P(1)$$
$$= \sqrt{n}\{1 - \hat{F}(2a - x) - \theta_0\} - 2f(2a - x)\sqrt{n}(\hat{a} - a) + o_P(1)$$
$$= \sqrt{n}\{1 - \hat{F}(2a - x) - \theta_0\} - 2f(x)f(a)^{-1}\sqrt{n}\{\hat{F}(a) - 1/2\} + o_P(1)$$

provided $f(a) > 0$, and hence $n^{-1/2}\sum_{i=1}^{n} m(X_i, \theta_0, \hat{a})$ is asymptotically normal from the Cramér–Wold device and the central limit theorem. It is easily seen that the asymptotic variance matrix $V_1$ is given by

$$V_1 = \begin{pmatrix} \theta_0(1 - \theta_0), & -\eta_0^2 - f(x)f(a)^{-1}\eta_0 \\ -\eta_0^2 - f(x)f(a)^{-1}\eta_0, & \theta_0(1 - \theta_0) + f(x)^2 f(a)^{-2} + 2f(x)f(a)^{-1}\eta_0 \end{pmatrix}.$$

The elements of this matrix can be estimated by replacing $\theta_0$ by $\hat{F}(x)$ and plugging in kernel estimators for $f(x)$ and $f(a)$.

Finally, we check condition (A0) when $0 < \theta_0 < 1/2$; the case $1/2 < \theta_0 < 1$ is similar. We need to show that $P\{(0, 0)^{\mathrm{t}} \in \mathcal{C}_n\} \to 1$. First, $P\{\hat{a} > x\} \to 1$ so we can condition on the event that $\hat{a} > x$. Next, note that $m(X, \theta_0, \hat{a})$ takes only three possible values:

$$\begin{pmatrix} 1 - \theta_0 \\ -\theta_0 \end{pmatrix}, \qquad \begin{pmatrix} -\theta_0 \\ 1 - \theta_0 \end{pmatrix} \quad \text{or} \quad \begin{pmatrix} -\theta_0 \\ -\theta_0 \end{pmatrix},$$

each with positive probability. It can be easily seen that the origin $(0, 0)^{\mathrm{t}}$ is contained in the interior of the convex hull of these three points, from which the assertion follows.

3.2. *Integral of squared densities.* Let $X_1, \ldots, X_n$ be i.i.d. from an unknown density $f_0$ which is assumed to be uniformly continuous and nonuniform. The quantity $\theta_0 = \int f_0^2 \, dx$ is of interest for various problems related to nonparametric density estimation. The limit distribution of the Hodges–Lehmann estimator of location has variance proportional to $1/\theta_0^2$; see Lehmann (1983), page 383. Similarly, the power of the Wilcoxon rank test is essentially determined by the size of $\theta_0$; see Lehmann (1975), page 72.

Consider the estimating function $m(X, \theta, f) = f(X) - \theta$ and let $m_n = n^{-1/2} m$. As a plug-in for $f_0$, we employ a kernel density estimator $\widehat{f}(x) = n^{-1} \sum_{i=1}^n k_b(X_i - x)$, where $k_b(\cdot) = k(\cdot/b)/b$ is a scaled version of a symmetric and bounded kernel function $k$ using bandwidth $b = b_n$. [For discussion of methods for deciding on good kernel bandwidths, when the specific purpose is precise estimation of $\theta_0$, see Schweder (1975).] Define

$$V = \int (f_0 - \theta_0)^2 f_0 \, dx = \int f_0^3 \, dx - \left( \int f_0^2 \, dx \right)^2,$$

which is the asymptotic variance of $n^{-1/2} \sum_{i=1}^n m(X_i, \theta_0, f_0)$, and is positive since $f_0$ is nonuniform. We now show that (A2) holds with $V_2 = V$. Write

$$n^{-1} \sum_{i=1}^n m^2(X_i, \theta_0, \widehat{f}) = n^{-1} \sum_{i=1}^n \{\widehat{f}(X_i) - \theta_0\}^2 = \int \widehat{f}^2 \, d\widehat{F} - 2\theta_0 \widehat{\theta} + \theta_0^2,$$

in terms of the empirical distribution function $\widehat{F}$ and $\widehat{\theta} = n^{-1} \sum_{i=1}^n \widehat{f}(X_i) = \int \widehat{f} \, d\widehat{F}$. Then $\int \widehat{f} \, d\widehat{F}$ and $\int \widehat{f}^2 \, d\widehat{F}$ have the required limits in probability, $\int f_0^2 \, dx$ and $\int f_0^3 \, dx$, respectively, provided $b \to 0$ and $nb \to \infty$. This verifies (A2).

Checking (A1) requires a more precise study of

$$\widehat{\theta} = n^{-1} \sum_{i=1}^n \widehat{f}(X_i) = n^{-2} \sum_{i,j} k_b(X_i - X_j) = \frac{k(0)}{nb} + \frac{n-1}{n} \widehat{g}.$$

Here $\widehat{g} = \widehat{g}(0)$, where $\widehat{g}(y) = \binom{n}{2}^{-1} \sum_{i<j} \bar{k}_b(Y_{i,j}, y)$ is a natural kernel estimator of the density $g(y) = \int f(y+x) f(x) \, dx$ of the difference $Y_{i,j} = X_i - X_j$, and $\bar{k}_b(Y_{i,j}, y) = \frac{1}{2} \{k_b(Y_{i,j} - y) + k_b(Y_{i,j} + y)\}$. Hjort (1999), Section 7, shows that $\widehat{g}(y)$ has mean value $g(y) + \frac{1}{2} b^2 g''(y) \int u^2 k(u) \, du + o(b^2)$, with variance $(4/n)\{g^*(y) - g(y)^2\}$ plus smaller order terms, where $g^*(y) = (1/4)\{\bar{g}(y, y) + \bar{g}(y, -y) + \bar{g}(-y, y) + \bar{g}(-y, -y)\}$ and $\bar{g}(y_1, y_2)$ is the joint density of two related differences $(X_2 - X_1, X_3 - X_1)$. It follows that

$$n^{-1/2} \sum_{i=1}^n m(X_i, \theta_0, \widehat{f}) = \sqrt{n}(\widehat{\theta} - \theta_0)$$

has mean of order $O(1/(\sqrt{n} b) + \sqrt{n} b^2)$ and variance going to $4V$. This, in conjunction with the asymptotic theory of U-statistics, verifies (A1) with $U \sim N(0, 4V)$, under the conditions $\sqrt{n}b \to \infty$ and $\sqrt{n}b^2 \to 0$. (If $b = b_0 n^{-\alpha}$, we need $\frac{1}{4} < \alpha < \frac{1}{2}$.) For (A3), note that $\widehat{f}(x) \le b^{-1}k_{\max}$ for all $x$, where $k_{\max}$ is the maximum of $k(u)$. Hence $\max_{i \le n} |\widehat{f}(X_i) - \theta_0|$ is bounded by $b^{-1}k_{\max} + \theta_0$, which implies (A3), provided only that $\sqrt{n}b \to \infty$.

Finally, for (A0) we need to show that

$$P\left\{ \min_{1 \le i \le n} m(X_i, \theta_0, \widehat{f}) < 0 < \max_{1 \le i \le n} m(X_i, \theta_0, \widehat{f}) \right\} \to 1.$$

First, consider

$$\max_{1 \le i \le n} m(X_i, \theta_0, \widehat{f}) \ge \max_{1 \le i \le n} f_0(X_i) - \max_{1 \le i \le n} |\widehat{f}(X_i) - f_0(X_i)| - \theta_0.$$

Note that $\max_{1 \le i \le n} |\widehat{f}(X_i) - f_0(X_i)| \to 0$ a.s. by the uniform consistency of $\widehat{f}$, which holds for $b$ as above (and suitable kernels $k$) by Theorem A of Silverman (1978), where we have used the assumption that $f_0$ is uniformly continuous. An example of a suitable kernel is the standard normal density function. Also, $\max_{1 \le i \le n} f_0(X_i) \to_{\text{a.s.}} \sup_t f_0(t) > \theta_0$, since $f_0$ is continuous and nonuniform, so $P\{\max_{1 \le i \le n} m(X_i, \theta_0, \widehat{f}) > 0\} \to 1$. In a similar way we can consider $\min_{1 \le i \le n} m(X_i, \theta_0, \widehat{f})$. We may now conclude that $-2\log \mathrm{EL}_n(\theta_0, \widehat{f}) \to_d 4\chi_1^2$.

3.3. *Functionals of survival distributions.* Wang and Jing (2001) (henceforth WJ) developed a plug-in version of EL for a class of functionals of a survival function (including its mean) in the presence of censoring. Denote the survival and censoring distribution functions by $F$ and $G$, respectively. The parameter of interest is a linear functional of $F$ of the form $\theta = \theta(F) = \int_0^\infty \xi(t) \, dF(t)$, where $\xi(t)$ is a (known) nonnegative measurable function and $\theta(F)$ is assumed finite. The estimating function is $m_n = n^{-1/2}m$, with

$$m(Z, \Delta, \theta, G) = \frac{\xi(Z)\Delta}{1 - G(Z)} - \theta,$$

$Z = \min(X, Y)$, $\Delta = I\{X < Y\}$, $Y \sim G$. Here $X \sim F$ and $Y \sim G$ are assumed to be independent. The Kaplan–Meier estimator $\widehat{G}_n$ of the censoring distribution function $G$ plays the role of the plug-in estimator. The resulting estimator $\widehat{\theta}$ of $\theta_0$ takes the form of an inverse-probability-weighted average. Equivalently, $\widehat{\theta} = \theta(\widehat{F}_n)$, where $\widehat{F}_n$ is the Kaplan–Meier estimator of $F$; see Satten and Datta (2001) for further discussion and references.

The conditions (A0)–(A3) needed to apply Theorem 2.1 are now checked by referring to various parts of WJ's proof of their Theorem 2.1, the conditions of which we assume implicitly. For (A0) we need to make the further mild assumption that the distribution of $\xi(X)$ is nondegenerate (i.e., not concentrated at its mean $\theta_0$). Then,

$$\max_{1 \leq i \leq n} m(Z_i, \Delta_i, \theta_0, \widehat{G}_n) \geq \max_{1 \leq i \leq n} \xi(Z_i)\Delta_i - \theta_0,$$

which is strictly positive for $n$ sufficiently large a.s. Also,

$$\min_{1 \leq i \leq n} m(Z_i, \Delta_i, \theta_0, \widehat{G}_n) = -\theta_0 < 0$$

for $n$ sufficiently large a.s. This, together with the lower bound for the maximum, entails (A0). Condition (A1) is immediate from the lemma on page 524 of WJ, with $U \sim \mathrm{N}(0, V_1)$ and $V_1$ being the asymptotic variance of $\widehat{\theta}$. Condition (A2) is checked using a Glivenko–Cantelli argument almost identical to that used below for estimation of $V_2$, where $V_2 = \mathrm{E}m^2(Z, \Delta, \theta_0, G) < \infty$ by condition (C3) of WJ. Condition (A3) is the display immediately before (4.5) in WJ.

It remains to provide consistent estimators of $V_1$ and $V_2$, and we do this along the lines of Remark 2.2. Stute's (1996) jackknife estimator can be used for $\widehat{V}_1$. Under conditions (A4)–(A5), we have that $\widehat{V}_2 = n^{-1} \sum_{i=1}^{n} m^2(Z_i, \Delta_i, \widehat{\theta}, \widehat{G}_n)$ consistently estimates $V_2$, where we also use the consistency of $\widehat{\theta}$. To check (A4), assume that $G(\tau_H-) < 1$, where $\tau_H = \inf\{t : H(t) = 1\}$, and $H$ is the distribution function of $Z$. Choose a constant $c$ such that $G(\tau_H-) < c < 1$. Specify $\bar{\mathcal{H}}$ as the class of increasing nonnegative functions $h$ such that $h(\tau_H-) < c$ and $h(t) = h(\tau_H)$ for $t \geq \tau_H$. Now, $\sup_{0 \leq t < \tau_H} |\widehat{G}_n(t) - G(t)|$ is bounded by

$$\sup_{0 \leq t < \tau_H} |\widehat{G}_n(t) - G(t \wedge Z_{(n)})| + \sup_{0 \leq t < \tau_H} |G(t \wedge Z_{(n)}) - G(t)|$$

$$= \sup_{0 \leq t \leq Z_{(n)}} |\widehat{G}_n(t) - G(t)| + \sup_{Z_{(n)} < t < \tau_H} |G(Z_{(n)}) - G(t)| \to_{\mathrm{pr}} 0,$$

by uniform consistency of $\widehat{G}_n$ on the interval $[0, Z_{(n)}]$; see Wang (1987). Thus $P\{\widehat{G}_n \in \bar{\mathcal{H}}\} = P\{\widehat{G}_n(\tau_H-) < c\} \to 1$. The class $\{1/(1 - h) : h \in \bar{\mathcal{H}}\}$ is contained in the class of all monotone functions into $[0, 1/(1 - c)]$, which is Glivenko–Cantelli; see van der Vaart and Wellner (1996), page 149. Thus, using the preservation property of Glivenko–Cantelli classes under a continuous function [see van der Vaart and Wellner (2000)], it follows that $\mathcal{F}$, defined right after conditions (A4) and (A5), is Glivenko–Cantelli. Condition (A5) follows by noting that $\mathrm{E}|m^2(Z, \theta, h) - m^2(Z, \theta_0, h)|$ is bounded above by

$$\mathrm{E}(|m(Z, \theta, h) - m(Z, \theta_0, h)||m(Z, \theta, h) + m(Z, \theta_0, h)|)$$

$$\leq \|\theta - \theta_0\|\{\|\theta + \theta_0\| + 2\,\mathrm{E}|\xi(Z)|/(1 - c)\}$$

for $h \in \bar{\mathcal{H}}$.

3.4. *Error distributions in nonparametric regression.* Consider the model $Y = \mu(X) + \varepsilon$, where $X$ and $\varepsilon$ are independent, $\varepsilon$ has unknown distribution function $F_\varepsilon$, and $\mu(\cdot)$ is an unknown regression function. We now use our approach with bootstrap calibration to construct an EL confidence interval for $\theta_0 = F_\varepsilon(z) \in (0, 1)$, at a fixed point $z$. The same assumptions as in Akritas and Van Keilegom (2001) are imposed. In particular, $F_\varepsilon$ is assumed to be continuous, $\mu(\cdot)$ is smooth and $X$ is bounded. For simplicity we restrict $X$ to $(0, 1)$.

Consider the Nadaraya–Watson estimator $\widehat{\mu}(x) = \sum_{i=1}^{n} W_{n,i}(x; b_n)Y_i$, with weights $W_{n,i}(x; b_n) = k_{b,x}(X_i) / \sum_{j=1}^{n} k_{b,x}(X_j)$ in terms of a kernel function $k$ and scaled versions $k_{b,x}(u) = b^{-1}k((u-x)/b)$ thereof, with $b = b_n = b_0 n^{-2/7}$ a bandwidth sequence (other choices of the bandwidth are possible). The estimating function is $m_n = n^{-1/2}m$, where $m(X, Y, \theta, \mu) = I\{Y - \mu(X) \leq z\} - \theta$.

We now check the conditions of Theorem 2.1. First, (A1) follows from the asymptotic normality of $\widehat{\theta} = n^{-1} \sum_{i=1}^{n} I\{\widehat{\varepsilon}_i \leq z\}$ [with $\widehat{\varepsilon}_i = Y_i - \widehat{\mu}(X_i)$], given by Theorem 2 in Akritas and Van Keilegom (2001): $\sqrt{n}\{\widehat{F}_\varepsilon(z) - F_\varepsilon(z)\} = n^{-1/2} \sum_{i=1}^{n} m(X_i, Y_i, \theta_0, \widehat{\mu}) \to_d N(0, V_1)$ where $V_1$ is defined in their paper. Condition (A2) holds with $V_2 = \theta_0(1 - \theta_0)$, provided $0 < \theta_0 < 1$. Also, (A3) holds since the function $\sqrt{n}m_n$ is uniformly bounded by 1. Finally, (A0) is an immediate consequence of the fact that $P\{Y - \widehat{\mu}(X) \leq z\}$ (probability conditionally on the function $\widehat{\mu}$) converges to $F_\varepsilon(z)$, which follows from a Taylor expansion and the uniform consistency of $\widehat{\mu}$. Since $F_\varepsilon(z)$ is strictly between 0 and 1, it follows that

$$P\{\text{there exist } 1 \leq i, j \leq n \text{ such that } Y_i - \widehat{\mu}(X_i) \leq z \text{ and } Y_j - \widehat{\mu}(X_j) > z\} \to 1,$$

which yields (A0).

It remains to estimate $V_1$ and $V_2$. Note that $\widehat{V}_2 = \widehat{\theta}(1 - \widehat{\theta})$ consistently estimates $V_2$. However, $V_1$ is harder to estimate. A plug-in type estimator can be obtained by making use of the estimator of the error density in Van Keilegom and Veraverbeke (2002). Since this approach requires the selection of a new bandwidth, we prefer to use the bootstrap approach. We now check the conditions of Theorem 2.2. For (A4), set $\delta > 0$ and define

$$C^{1+\delta}(0, 1) = \{\text{differentiable } f : (0, 1) \to \mathbb{R}, \text{such that } \|f\|_{1+\delta} \leq 1\},$$

where

$$\|f\|_{1+\delta} = \max\{\|f\|_\infty, \|f'\|_\infty\} + \sup_{x,y} \frac{|f'(x) - f'(y)|}{|x - y|^\delta},$$

and $\|\cdot\|_\infty$ denotes the supremum norm. Careful examination of the proof of Lemma 1 in Akritas and Van Keilegom (2001) reveals that the class $\{I(\varepsilon \leq$ $z + f(X))\colon f \in C^{1+\delta}(0,1)\}$ is Donsker, which is, using the notation of that proof, equal to the class $\mathcal{F}_1$ with $d_2 \equiv 1$ and $z$ fixed. Therefore, also the class

$$\{I(\varepsilon \le z + f(X)) - \theta\colon f \in C^{1+\delta}(0,1), \theta \in [0,1]\}$$
$$= \{I\{Y - h(X) \le z\} - \theta\colon h \in \bar{\mathcal{H}}, \theta \in [0,1]\}$$

is Donsker, and hence Glivenko–Cantelli, where $\bar{\mathcal{H}} = \mathcal{H} = \mu + C^{1+\delta}(0,1)$, and $\bar{\mathcal{H}}$ is endowed with the supremum norm. As a consequence, the class $\mathcal{F}$, defined right after (A4) and (A5), is also Glivenko–Cantelli. Moreover, $P\{\hat{\mu} \in \bar{\mathcal{H}}\} \to 1$ by Propositions 3–5 in Akritas and Van Keilegom (2001). Condition (A5) is satisfied since for any $\delta_n \downarrow 0$,

$$\sup_{|\theta - \theta_0| \le \delta_n, h \in \bar{\mathcal{H}}} |\mathrm{E}m^2(X,Y,\theta,h) - \mathrm{E}m^2(X,Y,\theta_0,h)|$$

$$\le \delta_n \sup_{|\theta - \theta_0| \le \delta_n, h \in \bar{\mathcal{H}}} \mathrm{E}|2I\{Y - h(X) \le z\} - \theta - \theta_0| \to 0.$$

Next, let us calculate $\Gamma(\theta, h)[\bar{h} - h]$ for any $h, \bar{h} \in \mathcal{H}$. We find

$$\lim_{\tau \to 0} \{M(\theta, h + \tau(\bar{h} - h)) - M(\theta, h)\}/\tau$$

$$= \lim_{\tau \to 0} \tau^{-1} \int [F_{Y|x}(z + h(x) + \tau(\bar{h}(x) - h(x))) - F_{Y|x}(z + h(x))]\, dF_X(x)$$

$$= \int f_{Y|x}(z + h(x))(\bar{h}(x) - h(x))\, dF_X(x),$$

where $F_{Y|x}$ and $f_{Y|x}$ are the distribution and density function of $Y$ given $X = x$, and $F_X$ is the distribution function of $X$. Consequently,

$$\Phi_n = \sqrt{n}\left[ n^{-1} \sum_{i=1}^{n} I\{Y_i - \mu(X_i) \le z\} - \theta_0 \right.$$

$$\left. + n^{-1} \int f_{Y|x}(z + \mu(x)) \sum_{i=1}^{n} (k_{b,x}(X_i)Y_i - \mathrm{E}\{k_{b,x}(X)Y\})\, dx \right]$$

(3)
$$+ o_{\mathrm{pr}}(1)$$

$$= \sqrt{n}\left[ n^{-1} \sum_{i=1}^{n} I\{Y_i - \mu(X_i) \le z\} - \theta_0 \right]$$

$$+ \sqrt{n}\left[ n^{-1} \sum_{i=1}^{n} f_{Y|X_i}(z + \mu(X_i))Y_i - \mathrm{E}[f_{Y|X}(z + \mu(X))Y] \right]$$

$$+ o_{\mathrm{pr}}(1).$$

In a similar way, we obtain

$$\Phi_n^* = \sqrt{n}\left[ n^{-1} \sum_{i=1}^{n} I\{Y_i^* - \hat{\mu}(X_i^*) \le z\} - n^{-1} \sum_{i=1}^{n} I\{Y_i - \hat{\mu}(X_i) \le z\} \right]$$

$$(4) \qquad + \sqrt{n}\left[n^{-1}\sum_{i=1}^{n} f_{Y|X_i^*}(z + \widehat{\mu}(X_i^*))Y_i^* - \mathrm{E}^*[f_{Y|X^*}(z + \widehat{\mu}(X^*))Y^*]\right]$$

$$+ o_{P^*}(1).$$

Both (3) and (4) converge to zero-mean normal random variables [use, e.g., the Lindeberg condition to show the convergence of (4)]. We next show that the asymptotic variance of (4) converges in probability to the asymptotic variance of (3). To show this we restrict attention to the first term of (3) and (4) (the convergence of the variance of the second term and of the covariance between the two terms can be established in a similar way). Note that the variance of the first term of (3) respectively (4) equals $\theta_0(1 - \theta_0)$ respectively $n^{-1}\sum_{i=1}^{n} I\{Y_i - \widehat{\mu}(X_i) \le z\}[1 - n^{-1}\sum_{i=1}^{n} I\{Y_i - \widehat{\mu}(X_i) \le z\}]$. Since it follows from Lemma 1 in Akritas and Van Keilegom (2001) that

$$n^{-1}\sum_{i=1}^{n} I\{Y_i - \widehat{\mu}(X_i) \le z\} = \theta_0 + \sum_{i=1}^{n}[I\{Y_i - \mu(X_i) \le z\} - \theta_0]$$

$$+ P\{Y - \widehat{\mu}(X) \le z \mid \widehat{\mu}\} - \theta_0 + o_{\mathrm{pr}}(n^{-1/2})$$

$$= \theta_0 + o_{\mathrm{pr}}(1),$$

the result follows. Hence, (B1) is satisfied. For (B2) it suffices by Remark 2.8 to show that the class $\{I\{Y - h(X) \le z\} - \theta : 0 \le \theta \le 1, h \in \widetilde{\mathcal{H}}\}$ is Donsker, which we have already established before, and that

$$\mathrm{Var}[I\{Y - h(X) \le z\} - I\{Y - \mu(X) \le z\} - \theta + \theta_0]$$

is bounded by $K_1|\theta - \theta_0| + K_2\|h - \mu\|_\infty$ for some $K_1, K_2 \ge 0$. A similar derivation can be given for the bootstrap analogue of (B2). Next write

$$|M(\theta_0, \widehat{\mu}) - \Gamma(\theta_0, \mu)[\widehat{\mu} - \mu]|$$

$$= \left|P\{Y - \widehat{\mu}(X) \le z\} - \theta_0 - \int f_{Y|x}(z + \mu(x))\{\widehat{\mu}(x) - \mu(x)\}\,dF_X(x)\right|$$

$$= \left|\int [F_{Y|x}(z + \widehat{\mu}(x)) - F_{Y|x}(z + \mu(x))\right.$$

$$\left. - f_{Y|x}(z + \mu(x))\{\widehat{\mu}(x) - \mu(x)\}]\,dF_X(x)\right|$$

$$= \frac{1}{2}\left|\int f'_{Y|x}(z + \xi(x))\{\widehat{\mu}(x) - \mu(x)\}^2\,dF_X(x)\right| \le K\sup_x |\widehat{\mu}(x) - \mu(x)|^2,$$

for some $\xi(x)$ between $\mu(x)$ and $\widehat{\mu}(x)$, and for some positive $K$. This shows that (B3) holds. In a similar way, the bootstrap version of (B3) can be shown to hold. Finally, condition (B4) follows from, for example, Härdle, Janssen and Serfling (1988), and its bootstrap version can be established in a very similar way. It now follows that a $100(1 - \alpha)\%$ confidence interval for $F_\varepsilon(z)$ is given by $\{\theta : -2\log \mathrm{EL}_n(\theta, \widehat{\mu}) \geq e^*_{1-\alpha}\}$, where $e^*_{1-\alpha}$ is the $100(1-\alpha)\%$ percentile of the distribution of

$$n\left[n^{-1}\sum_{i=1}^n I\{Y_i^* - \widehat{\mu}^*(X_i^*) \leq z\} - \widehat{\theta}\right]^2 \Big/ \{\widehat{\theta}(1-\widehat{\theta})\}.$$

3.5. *Density estimation.* Let $X_1, \ldots, X_n$ be i.i.d. from an unknown density $f_0$, and suppose we are interested in estimating $\theta_0 = f_0(t)$, for $t$ fixed. We do this using the kernel density estimator $\widehat{f}_n(t) = n^{-1}\sum_{i=1}^n k_b(X_i - t)$, where $k_b(u) = b^{-1}k(b^{-1}u)$ is a $b$-scaled version of a symmetric, bounded kernel function $k$, supported on $[-1, 1]$. We choose here to employ bandwidths $b = b_n$ that satisfy $nb \to \infty$ and $nb^5 \to 0$. The rate $b = cn^{-1/5}$ (for some $c > 0$) is optimal for estimating $f_0(t)$, in the sense of minimizing the asymptotic mean squared error, but as we here aim at constructing confidence intervals, an undersmoothing rate is preferable. Hall and Owen (1993) constructed EL confidence bands for $f_0$, and Chen (1996) showed that the pointwise EL confidence intervals (with and without Bartlett correction) are more accurate than those based on the bootstrap.

Following these authors, we use the sequence of estimating functions $m_n(x, \theta) = n^{-1/2}b^{1/2}\{k_b(x-t) - \theta\}$, which do not involve plug-in. We now check the conditions of Theorem 2.1. For (A0), note that $\sqrt{n}b^{-1/2}\min_{1 \leq i \leq n} m_n(X_i, \theta_0) = -\theta_0 < 0$, and

$$\sqrt{n}b^{-1/2}\max_{1 \leq i \leq n} m_n(X_i, \theta_0) = \max_{1 \leq i \leq n}\frac{1}{b}k\left(\frac{X_i - t}{b}\right) - \theta_0 \to_{\text{a.s.}} \infty$$

provided $f_0$ is bounded away from 0 in a neighborhood of $t$. Condition (A1) can be checked under mild conditions on the density, as it follows from standard asymptotic theory for kernel density estimators that $\sum_{i=1}^n m_n(X_i, \theta_0) = (nb)^{1/2}\{\widehat{f}_n(t) - f_0(t)\}$ tends to N$(0, V_1)$, where

$$(5) \qquad V_1 = f_0(t)R(k) \quad \text{and} \quad R(k) = \int k(u)^2 \, du.$$

For (A2),

$$\sum_{i=1}^n m_n^2(X_i, \theta_0) = \frac{b}{n}\sum_{i=1}^n \{k_b(X_i - t) - \theta_0\}^2 = \frac{1}{nb}\sum_{i=1}^n k((X_i - t)/b)^2 + O_{\text{pr}}(b),$$

which converges to $f_0(t)R(k) = V_1$ in probability. For (A3), $\max_{i \leq n}|m_n(X_i, \theta_0)| = O((nb)^{-1/2}) = o(1)$, because $k$ is bounded and $nb \to \infty$.

3.6. *Survival function estimation for current status data.* Suppose there is a failure time of interest $T \sim F$, with survival function $S = 1 - F$ and density $f$, but we only get to observe $Z = (C, \Delta)$, where $\Delta = I\{T \leq C\}$ and $C \sim G$ is an independent check-up time (with density $g$). The observations are assumed to be i.i.d.

The nonparametric maximum likelihood estimator $S_n(t)$ of $S(t)$ exists. Groeneboom (1987) showed that $n^{1/3}\{S_n(t) - S(t)\}$ converges to a nonde-generate limit law. The limit is not distribution-free, however, and is un-suitable for providing a confidence region for $S(t)$. Banerjee and Wellner (2005) found a universal limit law for the likelihood ratio statistic, leading to tractable confidence intervals. Our approach based on estimating equations offers a simpler type of EL confidence region, and extends to the setting in which $T$ and $C$ are conditionally independent given a covariate (although for simplicity we restrict attention to the case of no covariates).

First consider estimation of a smooth functional of $S$ (such as its mean): $\theta_0 = \int_0^\infty k(u)S(u)\,du$, where $k\colon[0,\infty) \to \mathbb{R}$ is fixed. This parameter can be estimated at a $\sqrt{n}$-rate, $m_n(Z, \theta, F, g, k) = n^{-1/2}m(Z, \theta, F, g, k)$ is an efficient influence curve, where

$$m(Z, \theta, F, k) = \frac{k(C)(1 - \Delta)}{g(C)} - \theta$$
$$- \frac{k(C)\{1 - F(C)\}}{g(C)} + \int_0^\infty k(u)\{1 - F(u)\}\,du,$$

and, given suitable preliminary estimators $\widehat{F}$ and $\widehat{g}$ of $F$ and $g$, respectively, we have a plug-in estimating function $m(Z, \theta, \widehat{F}, \widehat{g}, k)$ that yields a consistent estimator of $\theta_0$ when either $\widehat{F}$ or $\widehat{g}$ is consistent; see van der Laan and Robins (1998).

Now consider estimation of $0 < \theta_0 = S(t) < 1$. Van der Vaart and van der Laan (2006) introduced a kernel-type estimator $S_{n,b}(t)$ and showed that $n^{1/3}\{S_{n,b}(t) - S(t)\} \to_d \mathrm{N}(0, V_1)$, for appropriate and positive $V_1$. Their approach is to replace $k$ above by $k_n = k_{b,t}$, a kernel function of bandwidth $b = b_n = b_0 n^{-1/3}$ centered at $t$. Here $k_{b,t}(u) = k((u - t)/b)/b$ in terms of a bounded density $k$ supported on $[-1, 1]$. This yields a sequence of (plug-in) estimating functions $m_n(Z, \theta, \widehat{F}, \widehat{g}) = n^{-2/3}m(Z, \theta, \widehat{F}, \widehat{g}, k_n)$, and the estima-tor is written as $S_{n,b}(t) = \mathbb{P}_n\psi(\widehat{F}, \widehat{g}, k_n)$, where $\mathbb{P}_n$ is the empirical measure, and $\psi(F, g, k_n)(Z) = m(Z, 0, F, g, k_n)$ is the influence curve. The asymptotic variance of $S_{n,b}(t)$ is $V_1 = b_0^{-1}\sigma^2 R(k)$, where $R(k)$ is as in (5) and $\sigma^2$ depends on $F$ and $g$, as well as on the limits of $\widehat{g}$ and $\widehat{F}$.

We adopt the same assumptions as van der Vaart and van der Laan. In particular, assume that $F$ is differentiable at $t$, and $g$ is twice continuously differentiable and bounded away from zero in a neighborhood of $t$. Also, $\widehat{g}$ and $\widehat{F}$ are assumed to belong to classes of functions having uniform entropy of order $(1/\epsilon)^V$, for some $V < 2$, with probability tending to 1, and $\widehat{g}$, or $\widehat{F}$, or both, are locally consistent at $t$.

Our result for estimating functions with plug-in gives

$$-2\log \mathrm{EL}_n(S(t), \widehat{F}, \widehat{g}, k_n) \to_d \chi_1^2.$$

Conditions (A0)–(A3) are easily checked by referring to van der Vaart and van der Laan's Theorem 2.1 and its proof. For (A0), note that

$$n^{2/3} m_n(Z_i, \theta_0, \widehat{F}, \widehat{g}) = \frac{k_n(C_i)}{\widehat{g}(C_i)}(\widehat{F}(C_i) - \Delta_i)$$
$$+ \left[\int_0^\infty k_n(u)\{1 - \widehat{F}(u)\}\, du - \theta_0\right].$$

The minimum and maximum over $i \leq n$ of the first term above tend a.s. to $-\infty$ and $+\infty$, respectively, since $0 < P(\Delta = 1) < 1$ and it is assumed that $\widehat{g}$ is bounded away from 0 in a neighborhood of $t$. The second term above stays bounded as $n$ tends to infinity, so (A0) holds. Next, note that $\sum_{i=1}^n m_n(Z_i, \theta_0, \widehat{F}, \widehat{g}) = n^{1/3}\{S_{n,b}(t) - S(t)\}$, so (A1) holds [with $V_1$ given by the asymptotic variance of $S_{n,b}(t)$]. For (A2), note that

$$\sum_{i=1}^n m_n^2(Z_i, \theta_0, \widehat{F}, \widehat{g}) = n^{-1/3}\mathbb{P}_n\{\psi(\widehat{F}, \widehat{g}, k_n) - S(t)\}^2$$

(6)
$$= n^{-1/3}\mathbb{P}_n\{\psi(\widehat{F}, \widehat{g}, k_n) - P\psi(\widehat{F}, \widehat{g}, k_n)\}^2$$
$$+ 2n^{-1/3}\{S_{n,b}(t) - S(t)\}\{P\psi(\widehat{F}, \widehat{g}, k_n) - S(t)\}$$
$$- n^{-1/3}\{P\psi(\widehat{F}, \widehat{g}, k_n) - S(t)\}^2.$$

The last two terms above are asymptotically negligible, by the usual argument for controlling the bias of a kernel estimator; see the start of the proof of Theorem 2.1 of van der Vaart and van der Laan. To handle the first term, the influence function $\psi$ is split into a sum of two terms $\psi_1$ and $\psi_2$, where

$$\psi_2(F, g, k_n)(Z) = \int_0^\infty k_n(u)\{1 - F(u)\}\, du$$

does not give any contribution in the limit. In our case, $\psi_2$ acts as a constant function (there are no covariates), so the first term in (6) with $\psi$ replaced by $\psi_2$ is $O(n^{-1/3})$. The first term of (6) with $\psi$ replaced by $\psi_1$ can be expressed as

(7)
$$n^{-1/3}(n^{-1/2}\mathbb{G}_n H_n) + n^{-1/3} P H_n,$$

where $\mathbb{G}_n = \sqrt{n}(\mathbb{P}_n - P)$ is the empirical process and

$$H_n(\widehat{F}, \widehat{g}, k_n)(\cdot) = \{\psi_1(\widehat{F}, \widehat{g}, k_n) - P\psi_1(\widehat{F}, \widehat{g}, k_n)\}^2.$$

Applying the part of their proof that deals with $\psi_1$, but with $\psi_1$ replaced by $H_n$ and $n^{-1/2}k_n^2$ as the envelope functions, shows that $n^{-1/2}\mathbb{G}_n H_n$ is asymptotically tight. They also show that $n^{-1/3} PH_n \to_{\mathrm{pr}} b_0^{-1} \sigma^2 R(k)$, with $R(k)$ as in (5). Thus, only the second term in (7) gives a contribution in the limit, and we have

$$\sum_{i=1}^n m_n^2(Z_i, \theta_0, \widehat{F}, \widehat{g}) \to_{\mathrm{pr}} b_0^{-1} \sigma^2 R(k) = V_1,$$

establishing (A2) with $V_2 = V_1$. Finally, (A3) is checked using the assumption that $\widehat{g}$ is bounded away from zero in a fixed neighborhood of $t$. Note that $k_n \le c b_n^{-1} 1_{[t-b_n, t+b_n]}$ for some constant $c$, so

$$\max_{1 \le i \le n} |m_n(Z_i, \theta_0, \widehat{F}, \widehat{g})| = O_{\mathrm{pr}}(n^{-1/3}) = o_{\mathrm{pr}}(1).$$

**4. Empirical likelihood asymptotics with growing dimensions.** The traditional empirical likelihood theory works for a fixed number of estimating functions $p$, or, when estimating a mean, for data having a fixed dimension $d$. The present section is concerned with the question of how this theory may be extended toward allowing $p$ to increase with growing sample size. Consider situations with, say, $d$-dimensional observations $Z_1, \ldots, Z_n$ for which there are $p$-dimensional estimating functions $m(Z_i, \theta)$ to help assess a $p$-dimensional parameter $\theta$, and define

$$(8) \quad \mathrm{EL}_n(\theta) = \max\left\{ \prod_{i=1}^n (n w_i) \colon \text{each } w_i > 0, \sum_{i=1}^n w_i = 1, \sum_{i=1}^n w_i m(Z_i, \theta) = 0 \right\}.$$

Thus the framework is "triangular," reflecting a setup where the key quantities $p = p_n$, $d = d_n$, $Z_i = Z_{n,i}$, $\theta = \theta_n$, $m(z, \theta) = m_n(z, \theta)$ depend on $n$, but where we most of the time do not insist on keeping the extra subscript in the notation. A particular example would be $p$-dimensional $Z_i$'s for which their mean parameter $\mu$ is to be assessed, corresponding to estimation equation $m(z, \mu) = z - \mu$. We allow $p$ to grow with $n$, and study the problem of establishing sufficient conditions under which the standard $\chi_p^2$ calibration can still be used. There would often be a connection between $d$ and $p$, and indeed sometimes $d = p$, but the main interplay is between $n$ and $p$, and we do not need to make explicit requirements on $d = d_n$ itself.

We shall use several steps to approximate the EL statistic (8), and approximation results will be reached under different sets of conditions. Our results and tools for proving them shall involve the quantities

$$(9) \quad \bar{X}_n = n^{-1} \sum_{i=1}^n X_{n,i}, \qquad S_n = n^{-1} \sum_{i=1}^n X_{n,i} X_{n,i}^{\mathrm{t}}, \qquad D_n = \max_{i \le n} \|X_{n,i}\|,$$

where $X_{n,i} = m(Z_{n,i}, \theta_n)$. Here $\theta_n$ is the correct parameter, assumed to be properly defined as a function of the underlying distribution of $Z_{n,1}, \ldots, Z_{n,n}$ and the requirement that the mean value of $n^{-1} \sum_{i=1}^{n} m_n(Z_{n,i}, \theta_n)$ is zero (stressing in our notation, for this occasion, the dependence on $n$). We need $S_n$ to be positive definite, that is, at least $p$ among the $n$ vectors $X_{n,i}$ are linearly independent. In particular, $n \geq p$, and $p$ shall in fact have to grow somewhat slowly with $n$ in order for our approximation theorems to hold.

4.1. *Main results.* At the heart of the standard large-sample EL theorem lies the fact that

$$(10) \qquad T_n = -2 \log \mathrm{EL}_n(\theta_n) \quad \text{is close to} \quad T_n^* = n \bar{X}_n^{\mathrm{t}} S_n^{-1} \bar{X}_n.$$

One may view (10) as half of the story of how the EL behaves for large $n$ and $p$, the other half being how close $T_n^*$ then is to a $\chi_p^2$. A natural aim is therefore to secure conditions under which

$$(11) \qquad (T_n - T_n^*)/p^{1/2} \to_{\mathrm{pr}} 0 \quad \text{and} \quad (T_n^* - p)/(2p)^{1/2} \to_d \mathrm{N}(0, 1).$$

These statements taken together of course imply $(T_n - p)/(2p)^{1/2} \to_d \mathrm{N}(0, 1)$. Even though $(T_n - p)/(2p)^{1/2} \to_d \mathrm{N}(0, 1)$ may be achieved without (11), in special situations, we consider the quadratic approximation part and parcel of the EL distribution theory, and find it natural here to take "EL works for large $n$ and $p$" to mean both parts of (11).

Various sets of conditions may now be put up to secure (11), depending on the nature of the $X_{n,i}$ of (9). The following result provides an easily stated sufficient condition for (11) in the i.i.d. case, and has a number of applications that will be discussed in the next section.

THEOREM 4.1. *Suppose that the $X_{n,i}$'s are i.i.d. with mean zero and variance matrix $\Sigma_n$. First, if all components of $X_{n,i}$ are uniformly bounded and the eigenvalues of $\Sigma_n$ stay away from zero and infinity, then $p^3/n \to 0$ implies (11). Second, in case the components are not bounded, assume they have a uniformly bounded $q$th moment, for some $q > 2$, and again that the eigenvalues of $\Sigma_n$ stay away from zero and infinity. Then $p^{3+6/(q-2)}/n \to 0$ implies (11).*

The complete proof of Theorem 4.1 involves separate efforts for the two parts of (11), each of interest in its own right. We first explain the main ingredients in what makes the first part go through.

Introduce the random concave functions

$$(12) \quad G_n(\lambda) = 2 \sum_{i=1}^{n} \log(1 + \lambda^{\mathrm{t}} X_{n,i}/\sqrt{n}) \quad \text{and} \quad G_n^*(\lambda) = 2\lambda^{\mathrm{t}} \sqrt{n} \bar{X}_n - \lambda^{\mathrm{t}} S_n \lambda.$$

These are similar to the two random functions worked with in Remark 2.7, but are here defined in a somewhat different context. It is to be noted that $T_n$ of (10) is the same as $\max G_n = G_n(\widehat{\lambda})$, say, where the maximizer $\widehat{\lambda}$ also is the solution to $\sum_{i=1}^{n} X_{n,i}/(1 + \lambda^{\mathrm{t}} X_{n,i}/\sqrt{n}) = 0$. On the other hand, the maximizer of $G_n^*$ is $\lambda^* = S_n^{-1}\sqrt{n}\bar{X}_n$, and its maximum is precisely $T_n^*$. While $G_n^*$ is defined over all of $\mathbb{R}^p$, a little care is required for $G_n$, which is defined only where $\lambda^{\mathrm{t}} X_{n,i}/\sqrt{n} > -1$ for $i = 1, \ldots, n$. In view of the $(p/\sqrt{n})D_n \to_{\mathrm{pr}} 0$ condition that we nearly always shall impose, the (12) formula for $G_n$ holds with probability going to 1 for all $\lambda$ of size $O(p)$. We now provide basic "generic form" conditions for the first part of (11) to hold:

(D0) $P\{\mathrm{EL}_n(\theta_n) = 0\} \to 0$.
(D1) $(p/\sqrt{n})D_n \to_{\mathrm{pr}} 0$.
(D2) $\|\widehat{\lambda}\| = O_{\mathrm{pr}}(p^{1/2})$.
(D3) $\|\lambda^*\| = O_{\mathrm{pr}}(p^{1/2})$.
(D4) $\max \mathrm{eig}(S_n) = O_{\mathrm{pr}}(1)$.

PROPOSITION 4.1. *Conditions* (D0)–(D4) *imply* $(T_n - T_n^*)/p^{1/2} \to_{\mathrm{pr}} 0$. *If in addition* $(p^{3/2}/\sqrt{n})D_n \to_{\mathrm{pr}} 0$ *in* (D1), *then* $T_n - T_n^* \to 0$. *Furthermore, for both situations dealt with in Theorem 4.1, the conditions given there imply* (D0)–(D4).

Let us next focus on the second part of (11). Assume there is a population version $\Sigma_n$ of $S_n$ and consider $T_n^0 = n\bar{X}_n^{\mathrm{t}}\Sigma_n^{-1}\bar{X}_n$; when the $X_{n,i}$ are i.i.d., then $\Sigma_n$ is their variance matrix. Define

$$(13) \qquad L_n = |S_n - \Sigma_n| = \max_{j,k} |S_{n,j,k} - \Sigma_{n,j,k}|.$$

When $L_n$ is small, a well-behaved $\Sigma_n$ leads to a well-behaved $S_n$. We note that for any unit vector $u$, $|u^{\mathrm{t}} S_n u - u^{\mathrm{t}} \Sigma_n u| \leq \sum_{j,k} |u_j u_k| L_n \leq p L_n$, implying in particular that the range of eigenvalues for $S_n$ is within $pL_n$ of the range of eigenvalues for $\Sigma_n$. Also, $\mathrm{Tr}(S_n)$ is within $pL_n$ of $\mathrm{Tr}(\Sigma_n)$. Now consider the following conditions:

(D5) $p^{3/2}L_n \to_{\mathrm{pr}} 0$.
(D6) The eigenvalues of $\Sigma_n$ stay away from zero and infinity.

PROPOSITION 4.2. *Conditions* (D5)–(D6) *imply* $(T_n^* - T_n^0)/p^{1/2} \to_{\mathrm{pr}} 0$. *Furthermore, the assumptions detailed in Theorem 4.1 imply* (D5)–(D6), *for each of the two situations. Also, in the i.i.d. case, provided* $\mathrm{E}|X_{n,i,j}|^6$ *stays bounded for all components* $j \leq p$, *then the weak condition* $p/n \to 0$ *secures approximate* $\chi_p^2$*-ness in the sense that* $(T_n^0 - p)/(2p)^{1/2} \to_d \mathrm{N}(0, 1)$.

While Theorem 4.1 and corollaries indirectly noted above are satisfactory for several classes of problems, there are other situations of interest where the smallest eigenvalues, of $\Sigma_n$ and $S_n$, go to zero. This will typically lead to condition (D3) failing. For this reason we provide a parallel result that demands less regarding the distribution of eigenvalues. For the case of i.i.d. variables $X_{n,i} = m_n(Z_i, \theta_n)$ of mean zero and variance matrix $\Sigma_n$, consider $X_{n,i}^* = \Sigma_n^{-1/2} X_{n,i}$, and let $S_n^*$ be the empirical variance matrix of these, that is, $S_n^* = n^{-1} \sum_{i=1}^n Z_i^* (Z_i^*)^{\mathrm{t}} = \Sigma_n^{-1/2} S_n \Sigma_n^{-1/2}$. The eigenvalues of $S_n^*$ are often more well-behaved than those of $S_n$.

PROPOSITION 4.3. *Consider the EL setup of (8), with $m(Z_i, \mu) = Z_i - \mu$, for inference about the mean $\mu_n$ of $Z_i$. The conclusions of Theorem 4.1 continue to hold, without the condition on eigenvalues for $\Sigma_n$, as long as the conditions there are met for the transformed variables $Z_{n,i}^* = \Sigma_n^{-1/2}(Z_{n,i} - \mu_n)$.*

For another remark of relevance, write $\gamma_{1,n}$ and $\gamma_{p,n}$ for the largest and smallest eigenvalues of $\Sigma_n$. Yet another version of our main result emerges by dividing the $Z_i$'s by $\gamma_{p,n}^{1/2}$, to avoid small eigenvalues. This gives a parallel result to those of Theorem 4.1 and Proposition 4.3, where the essential condition is that the ratio $\gamma_{1,n}/\gamma_{p,n}$ remains bounded. See in this connection also Owen [[2001], page 86] where stability of this ratio is crucial also for some problems associated with fixed $p$.

For the four applications given in Section 5, along with a broad variety of others, the above development suffices. There are nevertheless situations where further variations on the conditions are required. In the following subsection the requirements (D0)–(D6) are discussed and followed up with further conditions that suffice for the different requirements to hold. We also give some useful lemmas that partly are needed to prove Propositions 4.1 and 4.2, and hence the master Theorem 4.1, and partly give the opportunity to prove versions of (11) under sets of conditions outside those of i.i.d. structures, like in regression models.

4.2. *On verifying conditions* (D0)–(D6). The EL operation (8) degenerates if zero is outside the convex hull spanned by $X_{n,1}, \ldots, X_{n,n}$ in $\mathbb{R}^p$. This may happen more frequently in higher dimensions. Condition (D0) amounts to the EL giving a positive maximum, with probability tending to 1 with $n$, and we now discuss conditions that secure this. That zero is outside the convex hull corresponds to there being a unit vector $u$ for which $u^{\mathrm{t}} X_{n,i} > 0$ for each $i$. So zero is inside the interior of the convex hull if $H_n(u) < 0$ for each unit vector $u$, where $H_n(u) = \min_{i \leq n} u^{\mathrm{t}} X_{n,i}$. Thus condition (D0) is implied by

$$(14) \qquad P\left\{ \max_{u \in \mathcal{U}_p} H_n(u) < 0 \right\} \to 1 \qquad \text{as } n \to \infty,$$

where $\mathcal{U}_p$ is the set of unit vectors in $\mathbb{R}^p$. This and several later problems will be handled separately for two types of situations: (a) the components of $X_{n,i}$ remain uniformly bounded, and (b) the components may be unbounded, but reasonable moment conditions prevail. It will be useful to deal with (D1) in connection with (14), that is, (D0). Yet another useful regularity condition is as follows.

(D7) For some $q > 2$, the sequence of $\mathrm{E}\|X_{n,i}/p^{1/2}\|^q$ stays bounded; and for this $q$ it holds that $p^{3+6/(q-2)}/n \to 0$.

LEMMA 4.1.   (a) *If the components of $X_{n,i}$ remain uniformly bounded, then $p^3/n \to 0$ implies* (D1). (b) *If* (D7) *holds, then again* (D1) *holds.*

LEMMA 4.2.   *For the i.i.d. case, assume there exists a positive $\varepsilon$ such that $r_p(u, \varepsilon) = P\{u^t X_{n,i} > -\varepsilon\} \le r < 1$ for all $u \in \mathcal{U}_p$; in particular, this necessitates a positive lower bound for the eigenvalues of $\Sigma_n$. (a) If the components of $X_{n,i}$ are uniformly bounded, then the requirement $(p \log p)/n \to 0$ as $n \to \infty$ secures* (14), *that is,* (D0). (b) *Also* (D7) *implies* (14).

Next we assess the sizes of the maximizers $\widehat{\lambda}$ and $\lambda^*$ of $G_n$ and $G_n^*$. We also need to inspect the size of $L_n$ of (13).

LEMMA 4.3.   *Suppose that $\sqrt{n}\|\bar{X}_n\| = O_{\mathrm{pr}}(p^{1/2})$, that $\min \mathrm{eig}(S_n)$ stays away from zero in probability, and that* (D1) *holds. Then $\|\widehat{\lambda}\| = O_{\mathrm{pr}}(p^{1/2})$, that is,* (D2) *holds.*

Note for the i.i.d. case, where the $X_{n,i}$'s have a variance matrix $\Sigma_n$, then $n\|\bar{X}_n\|^2$ is of the required size $O_{\mathrm{pr}}(p)$ if only $\mathrm{Tr}(\Sigma_n/p)$ stays bounded.

LEMMA 4.4.   *For the i.i.d. case, assume that the $X_{n,i,j}$'s have finite $q$th-order moments, for some $q \ge 4$, and let $A_n(p, q) = p^{-1} \sum_{j=1}^{p} \mathrm{E}|X_{n,i,j}|^q$. Then, for a positive constant $c(q)$,*

$$P\{L_n \ge \varepsilon\} \le \frac{c(q)p^2}{\varepsilon^q n^{q/2}} A_n(p, q)^2 \qquad \text{for each positive } \varepsilon.$$

It follows that when $q$th moments are bounded, then $p^{2+4/q}/n \to 0$ secures $pL_n \to_{\mathrm{pr}} 0$ and in its turn well-behaved eigenvalues of $S_n$ under minimal conditions on those of $\Sigma_n$. Similarly, the $p^{3+4/q}/n \to 0$ condition ensures $p^{3/2} L_n \to_{\mathrm{pr}} 0$, that is, (D5). We note further that when the $X_{n,i}$'s are uniformly bounded, then $p^2/n \to 0$ implies $pL_n \to_{\mathrm{pr}} 0$, whereas $p^3/n \to 0$ implies $p^{3/2} L_n \to_{\mathrm{pr}} 0$. This may be shown using techniques of the proof of Lemma 4.4. In situations where the $X_{n,i}$'s have moments of all orders, the growth conditions here come close to those for the case of bounded variables. If they are normal, for example, then $\|X_{n,i}\|$ is bounded by a variable of the type $c(\chi_p^2)^{1/2}$ for a suitable $c$, and one may show that $p^3/n \to 0$ and $p^4/n \to 0$ again suffice for $T_n^* - T_n$ being respectively $o_{\mathrm{pr}}(p^{1/2})$ and $o_{\mathrm{pr}}(1)$.

LEMMA 4.5. *Suppose* $\|T_n^*\| = O_{\mathrm{pr}}(p)$ *and that* $\min \mathrm{eig}(S_n)$ *stays away from zero in probability. Then* $\|\lambda^*\| = O_{\mathrm{pr}}(p^{1/2})$, *that is,* (D3) *holds.*

For condition (D4) we note for the i.i.d. case that if $\max \mathrm{eig}(\Sigma_n)$ is bounded, and $pL_n \to_{\mathrm{pr}} 0$, then (D4) holds, in view of comments made after (13).

Verifying eigenvalue conditions, for either $S_n$ or $\Sigma_n$, is sometimes technically hard. A theorem of Bai and Yin (1993) works for the case of $Z_i$ having independent components with zero means and unit variances, in which case the linear growth condition $p/n \to y \in (0, 1)$ ensures that the smallest and largest of the eigenvalues of $S_n$ tend a.s. to $(1 - \sqrt{y})^2$ and $(1 + \sqrt{y})^2$, respectively. Inspection of their proof reveals that a version holds also when $y = 0$, namely that the smallest and largest of eigenvalues then tend in probability to 1. See also Bai (1999) and the ensuing discussion.

## 5. Applications with growing $p$.

This section provides some examples where there is a growing number of parameters, and where the theory developed in Section 4 guarantees that the empirical likelihood methodology still is applicable.

5.1. *Many independent means.* Suppose that $Z_1, \ldots, Z_n$ correspond to $p$ independent samples $Z_{1,j}, \ldots, Z_{n,j}$, with mean $\mu_{0,j}$ and standard deviation $\sigma_j$, for $j = 1, \ldots, p$, assuming for simplicity of presentation that the sample size is the same for each group $j = 1, \ldots, p$. EL may then be used to make simultaneous inference for the vector of mean parameters $\mu_0$. Consider the normalized random vector $U_i$ with components $(Z_{i,j} - \mu_{0,j})/\sigma_j$, which has mean zero and variance matrix $I_p$. Results of Section 4 imply that the EL works properly, even when $p$ grows, provided $p^3/n \to 0$ and that the $U_i$ components stay uniformly bounded, for example, via Proposition 4.3. This is secured by the eigenvalue distribution result of Bai and Yin (1993) mentioned above.

Similar results may be reached in other models with a growing number of mean type parameters. An example is analysis of variance with a large number of groups; cf. Akritas and Arnold (2000). Our theory also supports the use of EL theory when multiple comparisons between groups are made, since the variance matrix of a collection of such differences of means is well-behaved enough to have its eigenvalues away from zero and infinity; that is, Theorem 4.1 applies.

5.2. *Poisson regression.* Assume that $Y_i$ given $z_i$ is Poisson with parameter $\mu_i = \exp(z_i^{\mathrm{t}}\beta)$, where $z_i$ is a $p$-dimensional covariate vector and $\beta$ a parameter vector of the same length. EL may be used with $\mathrm{EL}_n(\beta)$ defined as in (8), via estimating equations $\sum_{i=1}^{n} w_i\{Y_i - \exp(z_i^{\mathrm{t}}\beta)\}z_i = 0$. Assume that the covariate vectors $z_i$ are i.i.d. from some distribution, which we for an easy concrete illustration take to be the standard $p$-variate normal, and let us postulate further that the sequence of $\beta$ vectors is such that $\|\beta\|^2 = \sum_{j=1}^{p} \beta_j^2$ remains bounded, as $n$ and $p$ are allowed to grow. This fits the setup of Section 4 with $X_{n,i} = (Y_i - \mu_i)z_i$, which have variance matrix $\Sigma_n = \exp(\frac{1}{2}\|\beta\|^2)(I_p + \beta\beta^{\mathrm{t}})$. We see from this that its eigenvalues all lie between $\exp(\frac{1}{2}\|\beta\|^2)$ and $\exp(\frac{1}{2}\|\beta\|^2)(1 + \|\beta\|^2)$. The conditions of Theorem 4.1 hold, for each even $q$, from which we conclude that (11) holds as long as $p^{3+\varepsilon}/n \to 0$, for some positive $\varepsilon$.

This example may be generalized in various ways. The only point about the $N_p(0, I_p)$ distributional assumption for the covariates here was to get an explicit and easy $\Sigma_n$ matrix, and variations are easily constructed. Second, results can be derived inside the more usual regression framework where $z_1, \ldots, z_n$ are considered known covariate vectors. Basically, this involves the variance matrices $\Sigma_n = n^{-1}\sum_{i=1}^{n} \mu_i z_i z_i^{\mathrm{t}}$ and $S_n = n^{-1}\sum_{i=1}^{n}(Y_i - \mu_i)^2 z_i z_i^{\mathrm{t}}$, in generalization of those worked with in Section 4. Under a Lindeberg condition, combined with the requirement that the $|z_i^{\mathrm{t}}\beta|$ are bounded uniformly as $p$ and $n$ grow (which means that all Poisson means should be bounded away from zero and infinity), one may prove that conditions (D0)–(D6) are fulfilled as long as $p^3/n \to 0$, using methods associated with proving Lemmas 4.1–4.5. Hence the desired conclusion (11) holds. Similar results are reached for other generalized linear regression setups.

5.3. *Testing $f = f_0$ via orthogonal expansions.* For i.i.d. data $X_1, \ldots, X_n$ from an unknown density, consider the growing class of models

$$f_p(x \mid a) = f_0(x)c_p(a_1, \ldots, a_p)^{-1}\exp\left\{\sum_{j=1}^{p} a_j\psi_j(x)\right\}.$$

Here $f_0$ is a "start density," around which one models a flexible log-linear structure for deviations, the $\psi_j$ functions are orthonormal w.r.t. $f_0$, that is, $\int f_0\psi_j\psi_k\,dx = \delta_{j,k}$, and $c_p$ is the appropriate normalizing constant. Here we can carry out EL analysis for $\xi = (\xi_1, \ldots, \xi_p)^{\mathrm{t}}$, where $\xi_j = \int f\psi_j\,dx$, and a growing $p$. This is done via the vectors $Z_i = (\psi_1(X_i), \ldots, \psi_p(X_i))^{\mathrm{t}}$. The eigenvalues of its variance matrix will typically be well behaved, with reasonable conditions on $f$, and there is stability of fourth-order moments if, for example, the $\psi_j$'s are bounded. Thus EL theory holds for analysis of the $\xi_j$'s, if $p^3/n \to 0$. Consider in particular the problem of testing $f = f_0$, which corresponds to the $a_j$'s being zero. The theory of Section 4 ensures that

$T_n = -2\log\mathrm{EL}_n(0) = 2\sum_{i=1}^{n}\log(1 + \widehat{\lambda}^{\mathrm{t}}Z_i/\sqrt{n})$ is at most $O_{\mathrm{pr}}(p^{1/2})$ away from $T_n^* = n\bar{\psi}_n^{\mathrm{t}}S_n^{-1}\bar{\psi}_n$, where $\bar{\psi}_n$ is the vector of averages $n^{-1}\sum_{i=1}^{n}\psi_j(X_i)$ and $S_n = n^{-1}\sum_{i=1}^{n}Z_iZ_i^{\mathrm{t}}$, provided only $p^3/n \to 0$, if the $\psi_j$'s are uniformly bounded. This is since the variance matrix under the null hypothesis is simply equal to $I_p$. Also, both $T_n$ and $T_n^*$ have null distributions close enough to a $\chi_p^2$, again by Theorem 4.1.

5.4. *Growing polynomial regression.* Consider the regression model

$$Y_i = \xi(X_i) + \varepsilon_i \qquad \text{for } i = 1, \ldots, n,$$

where the pairs $(X_i, \varepsilon_i)$ are i.i.d., with $X_i$'s coming from some density $f$ and the $\varepsilon_i$'s having mean zero and standard deviation $\sigma_0$. The main objective is to make inference about $\xi(x)$. We do not strive for the fullest generality in this application of our theory, and are content to work with the following scenario: $f$ is known (e.g., the uniform on the unit interval), and $\xi(x)$ may be expanded in terms of basis functions $\psi_0, \psi_1, \psi_2, \ldots$ that are orthonormal w.r.t. $f$, that is, $\int f\psi_j\psi_k\,dx = \delta_{j,k}$, and where we take $\psi_0 = 1$. We might for example take $\psi_j(x) = \phi_j(F(x))$ where the $\phi_j$'s are orthogonal w.r.t. the uniform on the unit interval and $F$ the c.d.f. of $f$. Hence $\xi(x) = \sum_{j=0}^{\infty} b_j\psi_j(x)$, where we assume that $\mathrm{E}\xi(X)^2 = \sum_{j=0}^{\infty} b_j^2$ is finite, and also that $\xi(x)$ is bounded.

In this setup, consider as $p$th-order model

$$Y_i = \xi_p(X_i) + \varepsilon_i' \qquad \text{with } \xi_p(x) = \sum_{j=0}^{p} b_j\psi_j(x) = (\psi^{(p)}(x))^{\mathrm{t}}b^{(p)},$$

where the residuals are $\varepsilon_i' = \sum_{j=p+1}^{\infty} b_j\psi_j(X_i) + \varepsilon_i$ with variance $\sigma_p^2 = \sigma_0^2 + \sum_{j=p+1}^{\infty} b_j^2$; including more terms in the regression structure makes the residuals smaller in size, and vice versa. Consider $Z_i = Y_i\psi^{(p)}(X_i)$, a vector of dimension $p+1$, with mean value seen to be $b^{(p)}$. We will consider conditions under which $-2\log\mathrm{EL}_n(b^{(p)})$, based on $Z_1, \ldots, Z_n$, can be approximated by a $\chi_{p+1}^2$ distribution.

The key to verifying the conditions of Theorem 4.1 lies in controlling the sizes of the eigenvalues of the variance matrix of $Z_i$, which may be written

$$\Sigma_n = \mathrm{E}Y_i^2\psi^{(p)}(X_i)\psi^{(p)}(X_i)^{\mathrm{t}} - b^{(p)}(b^{(p)})^{\mathrm{t}} = \sigma_0^2 I_p + \Omega_p,$$

where $I_p$ and $\Omega_p$ are of size $(p+1) \times (p+1)$ and where the elements of the nonnegative definite $\Omega_p$ matrix are $\int \xi(x)^2\psi_j(x)\psi_k(x)f(x)\,dx - b_jb_k$. The eigenvalues of $\Sigma_n$ take the form $\sigma_0^2 + \phi_j$, where the $\phi_j$'s are the eigenvalues of $\Omega_p$, and are hence bounded downward by $\sigma_0^2$. They are also bounded upward, since for any unit vector $u$, $u^{\mathrm{t}}\Omega_p u$ is bounded by $M^2\int(u_0\psi_0 + \cdots + u_p\psi_p)^2 f\,dx = M^2$, where $M$ bounds $|\xi(x)|$.

As explained in the Introduction, we may apply our results to produce confidence regions for a subset of the $b_j$ parameters, to test whether some of them are equal to zero, and to make inference for any linear combination. For example, suppose we are interested in a simultaneous confidence region for the regression function $\xi_p$ at certain locations $x_1, \ldots, x_q$. Even though $p$ may be very large, provided $q$ grows slowly enough with $n$ our results apply to the focus parameter $\phi = (\xi_p(x_1), \ldots, \xi_p(x_q))$, because $\phi = f(b^{(p)})$ is a linear function of $b^{(p)}$ and the confidence region can be based on the transformed data $f(Z_i) = (\sum_{j=0}^{p} Y_i \psi_j(X_i) \psi_j(x_l))_{l=1,\ldots,q}$. Focus parameters defined by nonlinear functions would also be of interest, but this is beyond the scope of the paper, even in the case of a one-dimensional parameter such as $\phi = \max_{l=1,\ldots,q} |\xi_p(x_l)|$.

## 6. Concluding remarks.

6.1. *Nonstandard limit distributions.* Here we give a toy example in which $T_n = -2 \log \mathrm{EL}_n$ has a limit distribution different from $U^t V_2^{-1} U$ in Theorem 2.1. Let $X_i \sim \mathrm{N}(\theta_0, \sigma_i^2)$ be independent, and suppose $\sum_{i=1}^{\infty} \sigma_i^2 < \infty$. Consider the unbiased estimating function $m(X, \theta) = X - \theta$. Using steps from the proof of our Theorem 2.1, it can be shown that $T_n \to_d T$, the maximum of the process $G(\lambda) = 2 \sum_{i=1}^{\infty} \log(1 + \lambda \sigma_i Z_i)$ over the random interval $|\lambda| < 1/D$, where $D = \max_{i \geq 1} \sigma_i |Z_i|$ and the $Z_i$ are independent standard normals. In this case, (A0)–(A2) hold [with a random limit in (A2)], but (A3), which is needed to dispose of the remainder term in the quadratic approximation to $T_n$, fails, hence the nonstandard limit.

6.2. *Weighted EL.* The basic EL setup can be generalized to allow for weights. In the framework of Section 2, we can place a weight $\tau_i$ in front of each term $m_n(X_i, \theta, \hat{h})$ in $\mathrm{EL}_n(\theta)$. This would be useful in situations where the $X_i$'s have different precision. Conditions sufficing for $-2 \log \mathrm{EL}_n(\theta_0)$ to converge in distribution are readily developed, paralleling (A0)–(A3).

6.3. *Joint convergence of maximum and maximizer.* Our proof of Theorem 2.1 (in the case $a_n = 1$) shows that $T_n = \sup_\lambda G_n(\lambda)$ with probability tending to 1, and $\hat{\lambda} = \mathrm{argmax}_\lambda G_n(\lambda) = O_{\mathrm{pr}}(1)$, where $G_n(\lambda) = 2 \sum_{i=1}^{n} \log(1 + \lambda^t X_{n,i})$. Appealing to Theorem 5.1 of Banerjee and McKeague (2007), we can then infer the more general result that $(\hat{\lambda}, T_n) \to_d (V_2^{-1} U, U^t V_2^{-1} U)$. On the computational side, the proof also indicates that maximization or equation-solving algorithms should work better with $\lambda^* = V_n^{-1} U_n$ as starting point, rather than, for example, zero.

## APPENDIX

Here we provide proofs of theorems and claims presented earlier in our article.

PROOF OF THEOREM 2.1. The basic steps and notation of the proof were given in Remark 2.7. It remains to show that $T_n - T_n^* = o_{pr}(a_n)$, where $T_n = \sup G_n$ and $T_n^* = \sup G_n^*$. First we determine the stochastic order of $\hat{\lambda}$. Write $\hat{\lambda} = \|\hat{\lambda}\| u$, in terms of a random unit vector $u$. As in Owen [2001, page 220] we have

$$(15) \qquad \|\hat{\lambda}\| (u^t V_n u - D_n u^t U_n) \leq u^t U_n,$$

where $D_n = \max_{i \leq n} \|X_{n,i}\|$. But $u^t V_n u \geq \min \text{eig}(V_n) = O_{pr}(a_n^{-1})$, $u^t U_n = O_{pr}(1)$ and $D_n u^t U_n = o_{pr}(a_n^{-1})$, so $\|\hat{\lambda}\| = O_{pr}(a_n)$. Moreover, $\lambda^* = V_n^{-1} U_n$ when $V_n$ is invertible, so $\lambda^*$ is of the same stochastic order $O_{pr}(a_n)$ as $\hat{\lambda}$.

Write $\log(1+x) = x - \frac{1}{2}x^2 + \frac{1}{3}x^3 h(x)$, with $|h(x)| \leq 2$ for $|x| \leq \frac{1}{2}$. This gives, for any $c > 0$ and $\|\lambda\| \leq c$,

$$(16) \qquad G_n(\lambda) = 2\lambda^t U_n - \lambda^t V_n \lambda + r_n(\lambda),$$

where

$$|r_n(\lambda)| \leq (2/3) \sum_{i=1}^{n} |(\lambda^t X_{n,i})^3| |h(\lambda^t X_{n,i})|$$

$$\leq (4/3)\|\lambda\| D_n \lambda^t V_n \lambda \leq (4/3)c^3 D_n \max \text{eig}(V_n),$$

provided $cD_n \leq \frac{1}{2}$. With $T_{n,c}$ and $T_{n,c}^*$ denoting the maxima of $G_n$ and $G_n^*$ over the ball $\Omega_n(c) = \{\lambda : \|\lambda\| \leq ca_n\}$, we have

$$|T_{n,c}/a_n - T_{n,c}^*/a_n| \leq (1/a_n) \max\{|r_n(\lambda)| : \|\lambda\| \leq ca_n\}$$

$$\leq (4/3)c^3 a_n D_n \max \text{eig}(a_n V_n),$$

as long as $ca_n D_n \leq \frac{1}{2}$. Choose $c$ big enough to have both $\hat{\lambda}$ and $\lambda^*$ inside $\Omega_n(c)$ with probability above $1 - \eta$, for some preassigned $\eta$. Then

$$P\{|T_n/a_n - T_n^*/a_n| \geq \varepsilon\}$$

$$\leq P\{(4/3)c^2 a_n D_n \max \text{eig}(a_n V_n) \geq \varepsilon\}$$

$$+ P\{\|\hat{\lambda}\| > ca_n\} + P\{\|\lambda^*\| > ca_n\} + P\{ca_n D_n > \tfrac{1}{2}\}.$$

Hence the lim-sup of the probability sequence on the left is bounded by $2\eta$. Since $\eta$ was arbitrary, $T_n/a_n$ and $T_n^*/a_n$ must have the same limit distribution, namely $U^t V_2^{-1} U$. □

PROOF OF THE CLAIM OF REMARK 2.2.    Conditions (A4) and (A5) with $a_n = 1$ imply that, given any real sequence $\delta_n \downarrow 0$,

$$\sup_{\|\theta - \theta_0\| \leq \delta_n, h \in \tilde{\mathcal{H}}} \left| \sum_{i=1}^{n} \{m_n^{\otimes 2}(X_i, \theta, h) - m_n^{\otimes 2}(X_i, \theta_0, h)\} \right| \to_{\mathrm{pr}} 0.$$

The consistency of $\hat{\theta}$ then implies

$$R_n = \sum_{i=1}^{n} \{m_n^{\otimes 2}(X_i, \hat{\theta}, \hat{h}) - m_n^{\otimes 2}(X_i, \theta_0, \hat{h})\} \to_{\mathrm{pr}} 0.$$

Thus

$$|\hat{V}_2 - V_2| \leq |R_n| + \left| \sum_{i=1}^{n} m_n^{\otimes 2}(\theta_0, \hat{h}) - V_2 \right| \to_{\mathrm{pr}} 0,$$

where we have used assumption (A2) for the last term, so $\hat{V}_2$ consistently estimates $V_2$.   □

PROOF OF THEOREM 2.2.    By (2), the singular value theorem applied to $V_2^{-1}$ and $\hat{V}_2^{-1}$, along with the Cramér–Wold theorem, it suffices to show that $\hat{V}_2 \to_{\mathrm{pr}} V_2$ and that

$$P^* \{\sqrt{n}[M_n^*(\hat{\theta}, \hat{h}^*) - M_n(\hat{\theta}, \hat{h})] \leq t\} - P\{U \leq t\} = o_{\mathrm{pr}}(1).$$

The former follows from Remark 2.2, under conditions (A4) and (A5). For the latter, define, for any sequences $\alpha_n^1, \alpha_n^2 \downarrow 0$,

$$A_{n,\alpha_n} = \left\{ |\hat{\theta} - \theta_0| \leq \alpha_n^1, \sup_t |B_n(t)| \leq \alpha_n^1, \sup \|C_n(\theta, h)\| \leq \alpha_n^2 n^{-1/2}, \right.$$
$$\left. \|\hat{h} - h_0\|_{\mathcal{H}} \leq \alpha_n^1 n^{-1/4} \right\},$$

where $B_n(t)$ respectively $C_n(\theta, h)$ is the expression between absolute values (norm-signs) in condition (B1) respectively (B2), and where the supremum for $C_n$ is taken over $\|\theta - \theta_0\| \leq \alpha_n^1, \|h - h_0\|_{\mathcal{H}} \leq \alpha_n^1$. Then, by conditions (B1), (B2), (B4) and the consistency of $\hat{\theta}$, $\alpha_n^1$ and $\alpha_n^2$ can be chosen such that $P(A_{n,\alpha_n}) \to 1$ as $n$ tends to infinity. Hence it suffices to establish the convergence in probability, conditionally on the event $A_{n,\alpha_n}$. It now follows from condition (B5) that

$$\|M_n^*(\hat{\theta}, \hat{h}^*) - M_n^*(\hat{\theta}, \hat{h}) - \Gamma(\hat{\theta}, \hat{h})[\hat{h}^* - \hat{h}]\|$$
$$= \|M_n(\hat{\theta}, \hat{h}^*) - M_n(\hat{\theta}, \hat{h}) - \Gamma(\hat{\theta}, \hat{h})[\hat{h}^* - \hat{h}]\| + o_{P^*}(n^{-1/2})$$
$$\leq c\|\hat{h}^* - \hat{h}\|_{\mathcal{H}}^2 + o_{P^*}(n^{-1/2}) = o_{P^*}(n^{-1/2}) \qquad \text{a.s.}$$

In a similar way it follows from (B2), (B3) and (B4) that

$$\|M_n(\theta_0, \widehat{h}) - M_n(\theta_0, h_0) - \Gamma(\theta_0, h_0)[\widehat{h} - h_0]\| = o_{\mathrm{pr}}(n^{-1/2}).$$

Hence condition (B1) implies that

$$\sqrt{n}\{M_n^*(\widehat{\theta}, \widehat{h}^*) - M_n(\widehat{\theta}, \widehat{h})\}$$
$$= \sqrt{n}\{M_n^*(\widehat{\theta}, \widehat{h}) - M_n(\widehat{\theta}, \widehat{h}) + \Gamma(\widehat{\theta}, \widehat{h})[\widehat{h}^* - \widehat{h}]\} + o_{P^*}(1) \qquad \text{a.s.}$$

has the same limiting distribution as

$$\sqrt{n}\{M_n(\theta_0, h_0) + \Gamma(\theta_0, h_0)[\widehat{h} - h_0]\} = \sqrt{n}M_n(\theta_0, \widehat{h}) + o_{\mathrm{pr}}(1),$$

which by condition (A1) converges to $U$. $\quad\square$

We note that Theorem 4.1 is an immediate consequence of Propositions 4.1 and 4.2. We now tend to proving these.

PROOF OF PROPOSITION 4.1. That the conditions of Theorem 4.1 secure conditions (D0)–(D4) follows from Lemmas 4.1–4.4, proven below. Here we show that these conditions imply $(T_n - T_n^*)/p^{1/2} \to_{\mathrm{pr}} 0$.

Using (12) and (16) we see that $G_n^*$ is the natural two-step Taylor expansion approximation of $G_n$, and that $G_n = G_n^* + r_n$ with

$$r_n(\lambda) = (2/3)\sum_{i=1}^n (\lambda^{\mathrm{t}} X_{n,i}/\sqrt{n})^3 h(\lambda^{\mathrm{t}} X_i/\sqrt{n}) \le (4/3)\|\lambda\|(D_n/\sqrt{n})\lambda^{\mathrm{t}} S_n \lambda$$

as long as $\|\lambda D_n/\sqrt{n}\| \le \frac{1}{2}$. Choose $c$ such that the set $\Omega_n(c) = \{\lambda : \|\lambda\| \le cp^{1/2}\}$ catches both $\widehat{\lambda}$ and $\lambda^*$, with probability at least $1 - \eta$ for all large $n$, where $\eta$ is any preassigned positive number. Then

$$|r_n(\lambda)| \le (4/3)c^3 p^{3/2} n^{-1/2} D_n \max \mathrm{eig}(S_n) \qquad \text{for all } \lambda \in \Omega_n(c),$$

with arguments similar to those used for proving Theorem 2.1. This implies

$$P\{|T_n - T_n^*|/p^{1/2} \ge \varepsilon\}$$
$$\le P\{(4/3)c^3 p n^{-1/2} D_n \max \mathrm{eig}(S_n) \ge \varepsilon\}$$
$$\quad + P\{cp^{1/2} n^{-1/2} D_n > \tfrac{1}{2}\} + P\{\widehat{\lambda} \notin \Omega_n(c)\} + P\{\lambda^* \notin \Omega_n(c)\}.$$

Accordingly, under (D1)–(D4), the lim-sup of the left-hand side sequence is bounded by $2\eta$, and is hence zero. The modified and stronger result $T_n - T_n^* \to_{\mathrm{pr}} 0$ follows similarly under the stronger assumption. $\quad\square$

PROOF OF PROPOSITION 4.2. That the conditions of Theorem 4.1 guarantee conditions (D5)–(D6) is a consequence of Lemma 4.4, proven below. Here we show that these imply $(T_n^* - T_n^0)/p^{1/2} \to_{\mathrm{pr}} 0$.

To this end, write $S_n = \Sigma_n + \varepsilon_n$, so that $S_n^{-1} \doteq \Sigma_n^{-1} - \Sigma_n^{-1} \varepsilon_n \Sigma_n^{-1}$ when the elements of $\Sigma_n^{-1} \varepsilon_n$ become uniformly small, which they do in view of (D5)–(D6). Hence

$$T_n^* - T_n^0 \doteq n \bar{X}_n^{\mathrm{t}} \Sigma_n^{-1} \varepsilon_n \Sigma_n^{-1} \bar{X}_n = W_n^{\mathrm{t}} E_n W_n,$$

where $W_n = \Sigma_n^{-1/2} \sqrt{n} \bar{X}_n$ is seen to have $\|W_n\| = O_{\mathrm{pr}}(p^{1/2})$ and $E_n = \Sigma_n^{-1/2} \varepsilon_n \Sigma_n^{-1/2}$ must have the property that $|u^{\mathrm{t}} E_n u| = O_{\mathrm{pr}}(pL_n)$ for each unit vector $u$. This proves the first claim. The second claim of the proposition follows, after a transformation to new variables $X_{n,i}' = \Sigma_n^{-1/2} m_n(Z_i, \theta_n)$ with mean zero and variance matrix the identity matrix $I_p$, from efforts of Portnoy (1988), who used a martingale central limit theorem.  □

PROOF OF LEMMA 4.1. When $|X_{n,i,j}| \leq M$ for all components, then $D_n \leq Mp^{1/2}$, proving part (a). For the general case, to gauge the size of $D_n$ we cannot appeal to arguments involving the Borel–Cantelli lemma, as Owen (2001), Chapter 11, could when analyzing the fixed $p$ situation. However, $P\{(p/\sqrt{n})D_n \geq \varepsilon\}$ is bounded by

$$\sum_{i=1}^n P\{\|X_{n,i}\| \geq \varepsilon\sqrt{n}/p\} \leq n \frac{p^{3q/2}}{n^{q/2}\varepsilon^q} \max_{i \leq n} \mathrm{E}\|X_{n,i}/p^{1/2}\|^q,$$

which is seen to imply (b) of the lemma.  □

PROOF OF LEMMA 4.2. Observe that $|H_n(u) - H_n(v)| \leq \|u - v\| D_n$. The full surface of the $p$-dimensional unit ball may be covered by the union of a finite number $C_{p,n}$ of rectangles with side length $\delta_n$, provided $C_{p,n}\delta_n^{p-1}$ is as big as $A_p = 2\pi^{p/2}/\Gamma(p/2)$, the surface area of the unit ball. Hence

$$\max_{u \in \mathcal{U}} H_n(u) \leq \max_{u \in \mathcal{U}_{p,n}} H_n(u) + \delta_n D_n = H_n^* + \delta_n D_n,$$

where $\mathcal{U}_{p,n}$ is the finite set in question. To show (14) we demonstrate

$$P\{H_n^* < -\varepsilon\} \to 1 \quad \text{and} \quad P\{\delta_n D_n \leq \varepsilon\} \to 1.$$

We need to choose $\delta_n$ so that the second requirement holds, and then check whether $P\{H_n^* \geq -\varepsilon\} \leq C_{p,n}r^n$ is sufficient to meet the first requirement. What is demanded is that $\log C_{p,n} + n \log r \to -\infty$, and this is seen to correspond to $\{p \log(1/\delta_n)\}/n \to 0$.

(a) For the bounded components part we have $D_n \leq Mp^{1/2}$ as with Lemma 4.1, and may take $\delta_n = \varepsilon/(Mp^{1/2})$. In this case, therefore, the $n^{-1}p \log p \to 0$ condition suffices for (14) to hold. (b) For this situation we take $\delta_n = p/\sqrt{n}$, guaranteeing by Lemma 4.1 that $P\{\delta_n D_n \leq \varepsilon\} \to 1$. Some analysis shows that $(p/n)\log(1/\delta_n) = n^{-1/2}x_n \log(1/x_n)$, with $x_n = p/\sqrt{n}$, which tends to zero.  □

PROOF OF LEMMA 4.3. Write $\widehat{\lambda} = \|\widehat{\lambda}\|u$, where the random $u$ has unit length. One may argue as in Owen (2001), Chapter 11.2, to reach

$$\|\widehat{\lambda}\|\{u^{\mathrm{t}}S_n u - (D_n/\sqrt{n})\sqrt{n}u^{\mathrm{t}}\bar{X}_n\} \leq \sqrt{n}u^{\mathrm{t}}\bar{X}_n.$$

Here there is a positive $\delta$ such that the event $u^{\mathrm{t}}S_n u \geq \delta$ has probability tending to 1, while $D_n u^{\mathrm{t}}\bar{X}_n \rightarrow_{\mathrm{pr}} 0$. The result follows.  □

PROOF OF LEMMA 4.4. For the components of the $p \times p$ matrix $\varepsilon_n = S_n - \Sigma_n$, a bounding operation gives

$$P\{|\varepsilon_{n,j,k}| \geq \varepsilon\} \leq \frac{\mathrm{E}|\sqrt{n}\varepsilon_{n,j,k}|^q}{(\sqrt{n}\varepsilon)^q} \leq \frac{c(q)v_{n,j,k}^{q/2}}{n^{q/2}\varepsilon^q},$$

for a constant $c(q)$, by results of von Bahr (1965). Here $v_{n,j,k} = \mathrm{E}(X_{n,i,j}X_{n,i,k})^2 - (\Sigma_{n,j,k})^2$ is the variance of $X_{n,i,j}X_{n,i,k}$. This may be further bounded by

$$v_{n,j,k} \leq (\mathrm{E}|X_{n,i,j}|^4)^{1/2}(\mathrm{E}|X_{n,i,k}|^4)^{1/2} \leq (\mathrm{E}|X_{n,i,j}|^q)^{2/q}(\mathrm{E}|X_{n,i,k}|^q)^{2/q}$$

for $q \geq 4$. This leads to

$$P\{L_n \geq \varepsilon\} \leq \sum_{j,k} c(q)\frac{\mathrm{E}|X_{n,i,j}|^q\mathrm{E}|X_{n,i,k}|^q}{n^{q/2}\varepsilon^q},$$

which is then seen to imply the lemma.  □

PROOF OF LEMMA 4.5. We work with the explicit expression for $\lambda^*$, which leads to a representation in the form of $S_n^{-1/2}W_n$, with $W_n = S_n^{-1/2}\sqrt{n}\bar{X}_n$. Here $\|W_n\|$ is precisely $(T_n^*)^{1/2}$, hence of size $O_{\mathrm{pr}}(p^{1/2})$, while $\|S_n^{-1/2}u\| = O_{\mathrm{pr}}(1)$ for all unit vector $u$. This proves the lemma.  □

PROOF OF PROPOSITION 4.3. The central point to note is that the empirical likelihood (8) is invariant with respect to the transformation that maps data $Z_i$ to $A_nZ_i$, where $A_n$ is any nonsingular nonrandom $p \times p$ matrix. If $\mathrm{EL}_n(A_n\mu \mid A_n)$ is the empirical likelihood computed on the basis of $Z_i' = A_nZ_i$, for the parameter $\tilde{\mu} = A_n\mu$, then $A_n$ cancels out of the defining equation $\sum_{i=1}^n w_i(A_nZ_i - A_n\mu) = 0$, showing that $\mathrm{EL}_n(\tilde{\mu} \mid A_n)$ is the same as $\mathrm{EL}_n(\mu)$ in (8), that is, independent of $A_n$ (and with the same maximizing $w_i$'s). The same is true for the quadratic approximation $T_n = n(\bar{Z}_n - \mu_n)^{\mathrm{t}}S_n^{-1}(\bar{Z}_n - \mu_n)$ of (10). We may in particular employ $A_n = \Sigma_n^{-1/2}$, where the resulting $A_nZ_i$ have variance matrix $I_p$. The proof of the lemma now follows using arguments similar to those needed for Theorem 4.1 but under the additional simplifying assumptions that $\Sigma_n = I_p$.  □

**Acknowledgments.** We are very grateful to two careful reviewers, their Associate Editor and the Editor, for constructive comments that spurred a better organization of our manuscript and a clearer presentation of our results. The three authors are also grateful to the biostatistical working group NorEvent at the University of Oslo and its director Professor Odd Aalen for support in connection with a Research Kitchen in Oslo, and to the Institut de Statistique at Louvain-la-Neuve for generously arranging joint research visits.

N. L. HJORT
DEPARTMENT OF MATHEMATICS
UNIVERSITY OF OSLO
P.B. 1053 BLINDERN N–0316
OSLO
NORWAY
E-MAIL: nils@math.uio.no

I. W. MCKEAGUE
DEPARTMENT OF BIOSTATISTICS
COLUMBIA UNIVERSITY
722 WEST 168TH STREET
NEW YORK, NEW YORK 10032
USA
E-MAIL: im2131@columbia.edu

I. VAN KEILEGOM
INSTITUTE OF STATISTICS
UNIVERSITÉ CATHOLIQUE DE LOUVAIN
VOIE DU ROMAN PAYS 20
B–1348 LOUVAIN-LA-NEUVE
BELGIUM
E-MAIL: vankeilegom@stat.ucl.ac.be